\documentclass[preprint,12pt]{elsarticle}

\usepackage{blindtext}
\usepackage{hyperref}
\usepackage{amssymb}
\usepackage{amsmath}

\RequirePackage{algorithm}
\RequirePackage[noend]{algpseudocode}
\usepackage{xcolor}
\usepackage{silence}
\begin{document}

\begin{frontmatter}



\title{Deep Policy Iteration  for High-Dimensional Mean Field Games }


\author[inst1]{Mouhcine Assouli\corref{cor1}}

\author[inst1]{Badr Missaoui}

\cortext[cor1]{Corresponding author.}

\affiliation[inst1]{organization={Moroccan Center for Game Theory, UM6P},
            addressline={Rocade}, 
            city={Rabat-Salé},
            postcode={11103}, 
            country={Morocco}}

\begin{abstract}
This paper introduces Deep Policy Iteration (DPI), a novel approach that integrates the strengths of Neural Networks with the stability and convergence advantages of Policy Iteration (PI) to address high-dimensional stochastic Mean Field Games (MFG). DPI overcomes the limitations of PI, which is constrained by the curse of dimensionality to low-dimensional problems, by iteratively training three neural networks to solve PI equations and satisfy forward-backwards conditions. Our findings indicate that DPI achieves comparable convergence levels to the Mean Field Deep Galerkin Method (MFDGM), with additional advantages. Furthermore, deep learning techniques show promise in handling separable Hamiltonian cases where PI alone is less effective. DPI effectively manages high-dimensional problems, extending the applicability of PI to both separable and non-separable Hamiltonians. 

\end{abstract}

\begin{keyword}
Mean Field Games \sep Deep Learning  \sep Policy Iteration  \sep  Non-Separable Hamiltonian 

\end{keyword}

\end{frontmatter}


\section{Introduction}

Mean Field Games theory, introduced by Lasry and Lions \cite{lasry2007mean}, provides a framework for analyzing Nash equilibria in differential games involving a large number of agents, similar to Evolutionary Game Theory \cite{perc2010coevolutionary}, which models the dynamic evolution of strategies within the population. This theory has garnered significant attention and has been extensively studied in various fields, such as autonomous vehicles \cite{huang2019game, shiri2019massive}, finance \cite{cardaliaguet2018mean, casgrain2019algorithmic}, anthropology \cite{capraro2024language}, economics \cite{achdou2017income, achdou2014partial, gomes2015economic}. MFGs is characterized by a mathematical formulation consisting of a system of coupled partial differential equations (PDEs). Specifically, the system comprises a forward-time Fokker-Planck equation (FP) and a backward-time Hamilton-Jacobi-Bellman equation (HJB), which govern the evolution of the population density $\rho$ and the value function $\phi$, respectively.  In the general case, the MFG system is described as
\pdfstringdefDisableCommands{\def\<def>-command{}}
\begin{equation}\label{1}
\left\{
\begin{array}{rrrrr}
-\partial_t \phi-\nu \Delta \phi + H(x,\rho,p) =0, \ in&E,\\
\partial_t\rho-\nu \Delta \rho - \operatorname{div} \left(\rho \nabla_p H(x,\rho,p) \right)=0, \ in&E, \\
\rho(0,x)=\rho_0(x), \ \ \phi(T,x)=g(x,\rho(T,x)), \ in&\Omega ,
\end{array}
\right.
\end{equation}
where, $E= [0,T] \times \Omega,$  $\Omega$ bounded subset of $\mathbb{R}^d$, $p = \nabla \phi$
 and $g$ denotes the terminal cost. The Hamiltonian $H$ with separable structure is defined as
\begin{equation}
H(x,\rho,p) = H_0(x,p)-f_0(x,\rho), 
\end{equation}
where $H_0$ is the infimum of the Lagrangian function, which is the Legendre transform of the Hamiltonian, minus the interaction function $f_0$ between the population of agents.

The solution of (\ref{1}) exists and is unique under the standard assumptions of convexity of $H$ in the second variable and monotonicity of $f$ and $g$ \cite{lin2020apac}, see also refs \cite{lasry2006jeux, lasry2007mean} for more details. For non-separable Hamiltonians, where the Hamiltonian of the MFG depends jointly on $\rho$ and $p$, the existence and uniqueness of the solution for MFGs of congestion type has been investigated by Achdou and Porretta in \cite{article} and Gomes et al. in \cite{gomes2015short}.

The numerical solution of (\ref{1}) holds significant importance in the practical application of MFG theory. However, due to the strong coupling and forward-backward structure of the two equations in (\ref{1}), they cannot be solved independently or jointly using simple forward-in-time methods. Extensive research has been conducted in this area, leading to the proposal of various methods with successful applications, as seen in \cite{achdou2020mean, lasry2007mean} and references therein. Additionally, classical techniques demonstrate elegant theoretical convergence in MFG with separable Hamiltonians, such as finite difference schemes or semi-Lagrangian schemes, which have been investigated in \cite{achdou2010mean, achdou2013mean} and \cite{carlini2014fully, carlini2014semi} respectively. Furthermore, fictitious play \cite{hadikhanloo2019finite, cardaliaguet2017learning, gianatti2023approximation}. Nevertheless, these methods often face challenges in terms of computational complexity, particularly when dealing with high-dimensional problems  \cite{hammer1962adaptive, doi:10.1126/science.153.3731.34}. 

Recently, Deep learning approaches have shown promising results in addressing the dimensionality challenge, as demonstrated by research such as \cite{doi:10.1073/pnas.1922204117, carmona2021deep}. Further advancements have been made in techniques utilizing Generative Adversarial Networks (GANs), as evidenced in \cite{lin2020apac, cao2020connecting}. Unfortunately, none of these methods adequately cover the case of a non-separable Hamiltonian. 

To the best of our knowledge, two promising numerical approaches have arisen for mean-field games featuring non-separable Hamiltonians. The first approach is the Policy Iteration (PI) method  \cite{lauriere2021policy}, which is a numerical method based on finite difference techniques but may not adequately address separable Hamiltonians (see Test 3 subsection \ref{AC}). The second approach is the MFDGM method \cite{ASSOULI2023113802}, which is a deep learning approach based on the Deep Galerkin Method. 

Following the work in \cite{cacace2021policy, lauriere2021policy}, we introduce the PI algorithms for (1) with periodic boundary conditions.
We first define the Lagrangian as the Legendre transform of $H$ :
$$
L(\rho, q)=\sup _{p \in \mathbb{R}^d}\{p \cdot q-H(\rho, p)\}
$$
\textbf{Policy Iteration Algorithm:}  Given $R>0$ and given a bounded, measurable vector field $q^{(0)}: [0, T] \times \mathbb{T}^d  \rightarrow \mathbb{R}^d$ with $\left|q^{(0)}\right| \leq R$ and $\left\|\operatorname{div} q^{(0)}\right\|_{L^r(E)} \leq R$, iterate:\\

\noindent (i) Solve
\begin{equation}\label{fp}   
\begin{cases}\partial_t \rho^{(n)}-\epsilon \Delta \rho^{(n)}-\operatorname{div}\left(\rho^{(n)} q^{(n)}\right)=0, & \text { in } E, \\ \rho^{(n)}(x, 0)=\rho_0(x) & \text { in } \Omega. \end{cases}
\end{equation}
(ii) Solve
\begin{equation}\label{hjb}
\begin{cases}-\partial_t \phi^{(n)}-\epsilon \Delta \phi^{(n)}+q^{(n)} D \phi^{(n)}- L\left(\rho^{(n)}, q^{(n)}\right)=0 & \text { in } E, \\ \phi^{(n)}(T, x)=\phi_T(x) & \text { in } \Omega, \end{cases}
\end{equation}

(iii) Update the policy
\begin{equation}\label{control}
q^{(n+1)}(t, x)=\arg \max _{|q| \leq R}\left\{q \cdot D \phi^{(n)}(t, x)-L\left(\rho^{(n)}, q\right)\right\} \quad \text { in } E.
\end{equation}

\emph{Contributions} 
This paper introduces a novel approach, Deep Policy Iteration (DPI), aimed at approximately solving high-dimensional stochastic Mean Field Games (MFG). Our contribution lies in combining two key elements: firstly, leveraging the stability and convergence benefits of the Policy Iteration (PI) algorithm, based on the Banach fixed point method, to ensure controlled convergence over iterations. Secondly, we harness recent advancements in machine learning by employing Neural networks to address the curse of dimensionality. Inspired by prior work \cite{sirignano2018dgm,raissi2018deep}, we utilize three neural networks iteratively trained to approximate the unknown solutions of equations (\ref{fp}-\ref{control}) while satisfying each equation and its associated forward-backward conditions. Although PI alleviates the limitation of non-separable Hamiltonians, a common issue in existing methods, its applicability remains restricted to low dimensions and may not adequately address separable Hamiltonians. Thus, DPI significantly extends the applicability of PI to high-dimensional cases for both separable and non-separable Hamiltonian. To evaluate the reliability and effectiveness of DPI, we conduct a series of numerical experiments comparing its performance with that of the MFDGM and PI methods.\\

\emph{Contents} The subsequent sections of this paper are organized as follows: Section 2 provides an introduction to our approach, outlining its key aspects. In Section 3, we conduct a thorough review of related existing methods. Moving on to Section 4, we explore the numerical performance of our proposed algorithms. To evaluate the efficacy of our method, we employ a straightforward analytical solution in Section 4.1. In Section 4.2, we put our method to the test using two high-dimensional examples. Furthermore, we employ a well-established one-dimensional traffic flow problem \cite{huang2019game}, renowned for its non-separable Hamiltonian, to compare DPI and PI in Section 4.3.

\section{Methodology}\label{Methodology}
We utilize the PI algorithm to simplify the complexity of the MFG system. Our proposed methodology involves two main steps:
In the first step, the algorithm computes the density and cost of a given policy. This is achieved by solving a set of coupled partial differential equations (PDEs) (\ref{fp}) and (\ref{hjb})  that describe the dynamics of the system. In the second step, the algorithm updates the policy established in the first step and computes a new policy that minimizes the expected cost (\ref{control}) given the density and value function. Here, we proceed by using a distinct neural network for each of the three variables: the density function, the value function, and the policy. 

To assess the accuracy of these approximations, we use a loss function based on the residual of each equation to update the parameters of the neural networks. These neural networks are trained based on the Deep Galerkin Method \cite{sirignano2018dgm}. Within each iteration, the DPI method repeats the two steps until convergence. At each iteration, the algorithm computes a better policy by taking into account the expected behavior of the other agents in the population. The resulting policy is a solution to the MFG that describes the optimal behavior of the agents in the population, given the interactions between them.\\

We initialize the neural networks as a solution to our system and define:
\begin{equation}
q_{\tau}(t,x)=N_{\tau}(t,x), \ \ \ \phi_{\omega}(t,x)=N_{\omega}(t,x),  \ \ \ \rho_{\theta}(t,x)=N_{\theta}(t,x).
\end{equation}
Our training strategy starts by solving (\ref{fp}) with a fixed policy. We compute the loss (\ref{los1}) at randomly sampled points $\{(t_{b},x_{b})\}_{b=1}^{B}$ from $E$, and  $\{x_{s}\}_{s=1}^{S}$ from $\Omega$.
\begin{equation}\label{los1}
    \textsc{Loss}_{total}^{(3)}=\textsc{Loss}^{(3)}+\textsc{Loss}_{cond}^{(3)},
\end{equation}
where
\begin{align*}
\textsc{Loss}^{(3)}&=\frac{1}{B} \sum_{b=1}^{B} \Big|\partial_t \rho_{\theta}(t_{b},x_{b})
     - \nu \Delta \rho_{\theta}(t_{b},x_{b})\\
     &\quad \quad \quad- \operatorname{div} \left(\rho_{\theta}(t_{b},x_{b})q_{\tau}(t_{b},x_{b})) \right)\Big|^2,
\end{align*}
and
\begin{equation*}
   \textsc{Loss}_{cond}^{(3)}=\frac{1}{S} \sum_{s=1}^{S}\Big|\rho_{\theta}(0,x_{s})-\rho_0(x_{s})\Big|^2. 
\end{equation*}

We then update the weights of $\rho_{\theta}$ by back-propagating the loss (\ref{los1}). We do the same to (\ref{hjb}) with the given $\rho_{\theta}$. We compute (\ref{los2}) at randomly sampled points $\{(t_{b},x_{b})\}_{b=1}^{B}$ from $E$, and  $\{x_{s}\}_{s=1}^{S}$ from $\Omega$,\\
\begin{equation}\label{los2}
\textsc{Loss}_{total}^{(4)}=\textsc{Loss}^{(4)}+\textsc{Loss}_{cond}^{(4)}, 
\end{equation}
where
\begin{align*}
 \textsc{Loss}^{(4)}&=\frac{1}{B} \sum_{b=1}^{B}\Big| \partial_t \phi_{\omega}(t_{b},x_{b})
     + \nu \Delta \phi_{\omega}(t_{b},x_{b})
     \\ &\quad \quad + q_{\tau}(t_{b},x_{b}))\nabla \phi_{\omega}(t_{b},x_{b}) - L(\rho_{\theta}(t_{b},x_{b}), q_{\tau}(t_{b},x_{b}))\Big|^2,  
\end{align*}
and
$$\textsc{Loss}_{cond}^{(4)}=\frac{1}{S} \sum_{s=1}^{S} \Big|\phi_{\omega}(T,x_{s})-g(x_{s},\rho_{\theta}(T,x_{s}))\Big|^2.$$
We then update the weights of $\phi_{\omega}$  by backpropagating the loss (\ref{los2}). Finally, we update $q_{\tau}$ by computing 
the loss (\ref{los3}) at randomly sampled points $\{(t_{b},x_{b})\}_{b=1}^{B}$ from $E$.
\begin{equation}\label{los3}
    \textsc{Loss}_{total}^{(5)}=\textsc{Loss}^{(5)}+\textsc{Loss}_{cond}^{(5)},
\end{equation}
where
\begin{align*}
\textsc{Loss}^{(5)}&=\frac{1}{B} \sum_{b=1}^{B}  \Big| L(\rho_{\theta}(t_{b},x_{b}), q_{\tau}(t_{b},x_{b})) - q_{\tau}(t_{b},x_{b}))\nabla \phi_{\omega}(t_{b},x_{b}) \Big|^2,
\end{align*}
and
\begin{equation*}
   \textsc{Loss}_{cond}^{(5)}=\frac{1}{B} \sum_{b=1}^{B} \Big|q_{\tau}(t_{b},x_{b}) - \nabla_pH(\rho_{\theta}(t_{b},x_{b}),\nabla \phi_{\omega}(t_{b},x_{b}))\Big|^2. 
\end{equation*}

To help readers better understand our methodology, we have summarized it in the form of an algorithm \ref{alg1}. For reproducibility, we provide detailed settings and a publicly available Python implementation at \url{https://github.com/Mouhcine-56/-Deep-Policy-Iteration}.

The theoretical convergence of DPI, relying on the convergence of the neural network approximation was previously analyzed in \cite{sirignano2018dgm, ASSOULI2023113802}, leveraging the universal approximation theorem. Additionally, the convergence of PI has been investigated \cite{cacace2021policy, lauriere2021policy}  using the Banach fixed point method. It is important to note that by experimenting with diverse network structures and training approaches, one can enhance the performance and robustness of the neural networks utilized in the model. Hence, selecting the best possible combination of architecture and hyperparameters for the neural networks is crucial for achieving the desired outcomes. Additionally, while the convergence of residual loss is monitored, it alone is insufficient to determine the convergence of the neural network. Hence, a key feature of DPI is its reliance on the fixed-point algorithm, which enables the control of neural network convergence by observing the difference between the previous and current solutions in each iteration.
\begin{algorithm}[hbt!]
\begin{algorithmic}
\Require $L$, $\nu $ diffusion parameter, $g$ terminal cost. 
\Require  Initialize neural networks $ N_{\tau_0}$, $ N_{\omega_0}$ and $N_{\theta_0}$\\
\textbf{Train}
\For{n=0,1,2...,K-1}\\
     \State (i) solve (\ref{fp})
     \State Sample batch $\{(t_{b},x_{b})\}_{b=1}^{B}$ from $E$, and  $\{x_{s}\}_{s=1}^{S}$ from $\Omega$
     \State $\textsc{L}^{(3)}\leftarrow\frac{1}{B} \sum_{b=1}^{B} \Big|\partial_t \rho_{\theta_n}(t_{b},x_{b})
     - \nu \Delta \rho_{\theta_n}(t_{b},x_{b})$
      \State $\quad \quad \quad \quad \quad \quad \quad- \operatorname{div} \left(\rho_{\theta_n}(t_{b},x_{b})q_{\tau_n}(t_{b},x_{b})) \right)\Big|^2$
     \State $\textsc{L}_{cond}^{(3)}\leftarrow \frac{1}{S_1} \sum_{s=1}^{S}\Big|\rho_{\theta_n}(0,x_{s})-\rho_0(x_{s})\Big|^2.$
     \State Backpropagate $\textsc{Loss}_{total}^{(3)}$  to  $\theta_{n+1}$ weights.
     \State \\
     \State (ii) solve (\ref{hjb})
     \State Sample batch $\{(t_{b},x_{b})\}_{b=1}^{B}$ from $E$, and  $\{x_{s}\}_{s=1}^{S}$ from $\Omega$.
     \State $\textsc{L}^{(4)}\leftarrow \frac{1}{B} \sum_{b=1}^{B}\Big| \partial_t \phi_{\omega_n}(t_{b},x_{b})
     + \nu \Delta \phi_{\omega_n}(t_{b},x_{b})$
     \State $\quad \quad \quad \quad + q_{\tau_n}(t_{b},x_{b}))\nabla \phi_{\omega_n}(t_{b},x_{b})- L(\rho_{\theta_{n+1}}(t_{b},x_{b}), q_{\tau_n}(t_{b},x_{b}))\Big|^2$ 
     \State $\textsc{L}_{cond}^{(4)}\leftarrow\frac{1}{S} \sum_{s=1}^{S} \Big|\phi_{\omega_n}(T,x_{s})-g(x_{s},\rho_{\theta_{n+1}}(T,x_{s}))\Big|^2.$
     \State Backpropagate  $\textsc{Loss}_{total}^{(4)}$ to $\omega_{n+1}$  weight.
     \State \\
     \State (i) solve (\ref{control})
     \State Sample batch $\{(t_{b},x_{b})\}_{b=1}^{B}$ from $E$.
     \State $\textsc{L}^{(5)}\leftarrow \frac{1}{B} \sum_{b=1}^{B}  \Big|L(\rho_{\theta_{n+1}}(t_{b},x_{b}), q_{\tau_n}(t_{b},x_{b})) -q_{\tau_n}(t_{b},x_{b}))\nabla \phi_{\omega_{n+1}}(t_{b},x_{b})\Big|^2,$ 
     \State $\textsc{L}_{cond}^{(5)}\leftarrow\frac{1}{B} \sum_{b=1}^{B} \Big|q_{\tau_n}(t_{b},x_{b}) - \nabla_pH(\rho_{\theta_{n+1}}(t_{b},x_{b}),\nabla \phi_{\omega_{n+1}}(t_{b},x_{b}))\Big|^2.$
     \State Backpropagate  $\textsc{Loss}_{total}^{(5)}$ to $\tau_{n+1}$  weight.
     \State \\
\EndFor 
\Return $\theta_{K}$, $\omega_{K}$  and $\tau_{K}$.
\end{algorithmic}
\caption{Deep Policy Iteration Method }\label{alg1}
\end{algorithm}
\section{Related Works}\label{Related Works}
In this section, we present a literature review that pertains to our study. Specifically, we concentrate on pertinent sources that provide insight into our research.

{\bf GANs:} Generative adversarial networks, introduced in 2014 \cite{goodfellow2020generative}, represent a class of machine learning models known for their success in image generation and data processing \cite{denton2015deep, reed2016generative, radford2015unsupervised}. Recently, there has been a growing interest in employing GANs for financial modeling as well \cite{wiese2019deep}. The GAN architecture comprises two neural networks: a generator and a discriminator, operating in opposition to each other to produce samples from a specific distribution. The overarching objective is to achieve equilibrium in solving the following optimization problem,
\begin{equation}\label{gan6}
\min\limits_{G}\max\limits_{D}\Big\{\mathbb{E}_{x \sim P_{data}(x)}[log(D(x)]\\
     + \mathbb{E}_{z \sim P_{g}(z)}[log(1-D(G(z))]\Big\},
\end{equation}
In this formulation, $P_{\text{data}}(x)$ represents the distribution of the original data, and $P_{g}(z)$ represents the distribution of noise data. The primary objective is to minimize the output of the generator ($G$) and simultaneously maximize the output of the discriminator ($D$). This is achieved by comparing the probability of the discriminator correctly identifying original data $P_{\text{data}}(x)$ with the probability of the discriminator incorrectly classifying generated data $G$ produced by the generator with noise data $P_{g}(z)$ as real ($1-D(G(z))$).
Essentially, the discriminator aims to distinguish between real and fake data accurately, while the generator endeavours to create synthetic data that can convincingly deceive the discriminator. This adversarial process results in the generator progressively improving its ability to generate realistic data.\\

{\bf   APAC-Net:} In the work by Lin et al. \cite{lin2020apac}, the authors propose a method called APAC-Net based on Generative Adversarial Networks (GANs) to address high-dimensional MFGs in stochastic scenarios. They leverage the Hopf formula in density space to reformulate MFGs as a saddle-point problem given by, 
\begin{equation}\label{eqGM4}
\begin{split}
\inf\limits_{\rho(x,t)}\sup\limits_{\phi(x,t)}&\Big\{\mathbb{E}_{z\sim P(z),t \sim Unif[0,T]} [\partial_t \phi(\rho(t,z),t)+ \nu \Delta \phi(\rho(t,z),t
      )\\  &\quad- H(\rho(t,z),\nabla \phi)]
      + \mathbb{E}_{z\sim P(z)} \phi(0,\rho(0,z)) - \mathbb{E}_{x \sim \rho_T} \phi(T,x) \Big \}, 
\end{split}
\end{equation}
Remarkably, a noteworthy link emerges between GANs and MFGs, as Equation (\ref{eqGM4}) facilitates a connection to the Kantorovich-Rubenstein dual formulation of Wasserstein GANs \cite{villani2021topics}, expressed as:
\begin{equation}\label{gan6bis}
\begin{array}{cc}
     \min\limits_{G}\max\limits_{D} \{\mathbb{E}_{x \sim P_{data}(x)}[(D(x)] - \mathbb{E}_{z \sim P_{g}(z)}[(D(G(z))]\},  \\
     s.t. \ \ ||\nabla D||\leq 1. 
\end{array}
\end{equation}
 this formulation allows us to sidestep the use of spatial grids or uniform sampling in high dimensions, addressing the curse of dimensionality. While \cite{doi:10.1073/pnas.1922204117} also avoids spatial grids in Mean Field Games (MFGs), it is important to note that their study is limited to the deterministic setting when the parameter $\nu = 0$.
 
Consequently, an algorithm analogous to GANs can be employed to address MFG problems. However, this requires the Hamiltonian H to be separable. In such cases, solving the MFG-LWR system (details in Section \ref{tf}) becomes unfeasible. Generally, MFGs with non-separable Hamiltonians, preventing their reformulation as in Equation (\ref{eqGM4}), pose significant challenges for resolution with GANs.
\\

{\bf MFDGM:} The methodology proposed by \cite{ASSOULI2023113802} involves the use of two neural networks to approximate the population distribution and the value function. The accuracy of these approximations is optimized by using a loss function based on the residual of the first equation (HJB) to update the parameters of the neural networks. The process is then repeated using the second equation (FP) and new parameters to further improve the accuracy of the approximations. In contrast, methods based on Generative Adversarial
Networks (GANs) \cite{cao2020connecting, lin2020apac} cannot solve MFGs with non-separable Hamiltonians. This novel methodology significantly improves the computational complexity of the method by utilizing neural networks trained simultaneously with a single hidden layer and optimizing the approximations through a loss function. This makes the method more efficient and scalable, enabling its application to a wider range of high-dimensional problems in general mean-field games, including separable or non-separable Hamiltonians and deterministic or stochastic models. Additionally, comparisons to previous methods demonstrate the efficiency of this approach.\

Our approach may offer similar capabilities to MFDGM while providing the advantage of controlling convergence by integrating PI\\

{\bf DGM-MFG:} In the study by Carmona et al. \cite{carmona2021deep}, Section 4 intricately details the adaptation of the Deep Galerkin Method (DGM) algorithm for the resolution of Mean Field Games (MFGs). This method proves effective across a diverse spectrum of partial differential equations (PDEs), notable for its independence from the specific structural details of the problem.  MFDGM shares a fundamental similarity with DGM-MFG in the utilization of neural networks for approximating unknown functions. Both methods also involve parameter adjustment to minimize a loss function derived from the PDE residual, a commonality observed in both \cite{carmona2021deep} and \cite{cao2020connecting}. However, MFDGM stands out in its training methodology. Instead of relying on the aggregation of Partial Differential Equation (PDE) residuals as a unified loss function, MFDGM defines distinct loss functions for each equation.
\\

{\bf Policy Iteration method:} In \cite{ lauriere2021policy}, the authors introduced the Policy iteration method, which proved to be the first successful approach for solving systems of mean-field game partial differential equations with non-separable Hamiltonians. The method involved proposing two algorithms based on policy iteration, which iteratively updated the population distribution, value function, and control. Since the control was fixed, these algorithms only required the solution of two decoupled, linear PDEs at each iteration, which reduced the complexity of the equations; refer to \cite{cacace2021policy}. 

 The authors presented a revised version of the policy iteration algorithm, which involves updating the control during each update of the population distribution and value function. Their novel approach strives to make significant strides in the field and deepen our understanding of the subject matter. However, due to the computational complexity of the method, it was limited to low-dimensional problems.

One major limitation of this method was that it may not work in the separable case. To address these issues, our study proposes a new approach using neural networks. We will use the same policy iteration method but instead of using the finite difference method to solve the two decoupled, linear PDEs, we will use neural networks inspired by DGM. This will allow us to overcome the computational challenges of the original method and extend its applicability to the separable case.
With this approach, we expect to significantly improve the computational complexity of the method, which will allow us to apply it to a wider range of high-dimensional problems. 

Overall, our study seeks to extend the applicability of the PI method to more complex and higher-dimensional problems by incorporating neural networks into its framework. We believe that this approach has the potential to make a significant contribution to the field of mean-field game theory and could have practical implications in fields such as finance, economics, and engineering.\\

\section{Numerical Experiments}\label{Numerical Experiments}
To assess the efficacy of the proposed algorithm (\ref{alg1}), we conducted experiments on various problems. Firstly, we utilized the example dataset provided in 
\cite{lin2020apac,ASSOULI2023113802}, which has an explicitly defined solution structure that facilitates straightforward numerical comparisons. Secondly, we extended the application of our method to more intricate problems involving high-dimensional cases, further substantiating its reliability. Finally, we evaluated the algorithm's performance on a well-studied one-dimensional traffic flow problem \cite{huang2019game}, which is known for its non-separable Hamiltonian.\\
To ascertain the reliability of our approach, we compared the performance of three different algorithms, namely PI, DPI, and MFDGM, on the same dataset.  This choice was motivated by the fact that previous methodologies from the literature had already undergone comparisons with MFDGM \cite{ASSOULI2023113802}. Through this evaluation, we aimed to determine the effectiveness of our proposed algorithm in comparison to existing state-of-the-art methods.

\subsection{Analytic Comparison}\label{AC}
We test the effectiveness of the DPI method, we compare its performance with a simple example of the analytic solution characterized by a separable Hamiltonian. In order to simplify the comparison process, we select the spatial domain as $\Omega=[-2,2]^d$ and set the final time as $T = 1.$ This allows us to easily analyze and compare the results obtained from both methods. For
\begin{equation}
\begin{array}{cc}
     H_0(x,p)=\frac{||p||^2}{2}-\beta \frac{||x||^2}{2}, \ \ \ f_0(x,\rho)=\gamma \ln(\rho),  \\
     g(x)=\alpha\frac{||x||^2}{2}-(\nu d \alpha +\gamma\frac{d}{2}\ln \frac{\alpha}{2\pi\nu}),
\end{array}
\end{equation}
and $\nu=\beta=1$,
where $$\alpha=\frac{-\gamma +\sqrt{\gamma^2+4\nu^2 \beta}}{2\nu}=\frac{-\gamma +\sqrt{\gamma^2+4}}{2}.$$
The corresponding MFG system is:
\\
\begin{equation}\label{test}
\left\{
\begin{array}{rrrrr}
-\partial_t \phi- \Delta \phi + \frac{||\nabla \phi||^2}{2}- \frac{||x||^2}{2} = \gamma \ln(\rho), \\
\partial_t\rho- \Delta \rho - \operatorname{div} \left(\rho \nabla \phi \right)=0,  \\
\rho(0,x)=(\frac{1}{2\pi})^{\frac{d}{2}}e^{-\frac{ ||x||^2}{2}}, \\ \phi(T,x)=\alpha\frac{||x||^2}{2}-(  \alpha d +\gamma\frac{d}{2}\ln \frac{\alpha}{2\pi}),
\end{array}
\right.
\vspace{-0.25cm}
\end{equation}
\\
and the explicit formula is given by 
\begin{equation}
\begin{array}{cc}
     \phi(t,x)=\alpha\frac{||x||^2}{2}-(  \alpha d +\gamma\frac{d}{2}\ln \frac{\alpha}{2\pi})t, \\
     \rho(t,x)=(\frac{1}{2\pi})^{\frac{d}{2}}e^{-\frac{ ||x||^2}{2}}.
\end{array}
\end{equation}

\begin{figure}
\centering\includegraphics[width=3.8cm]{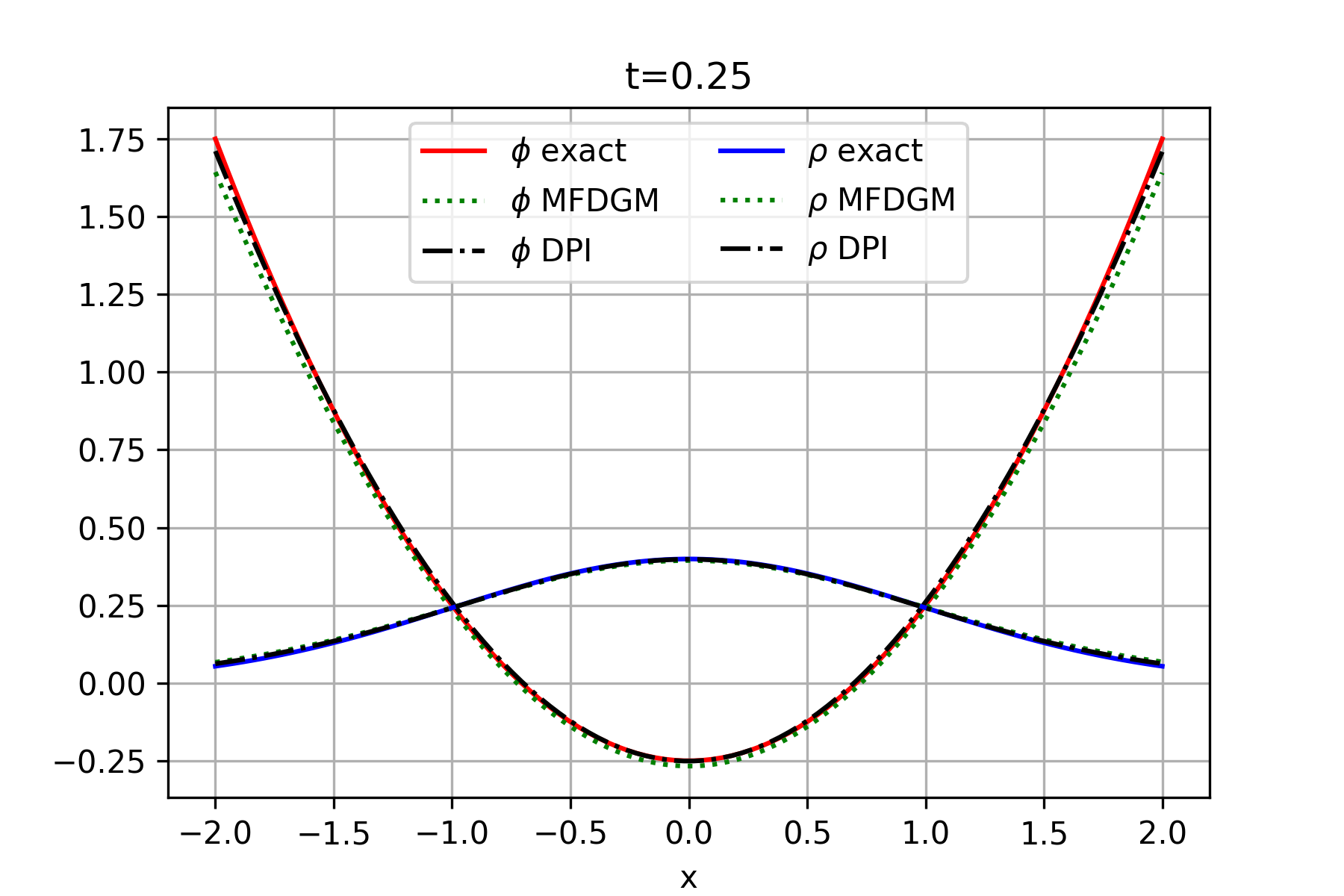}
\centering\includegraphics[width=3.8cm]{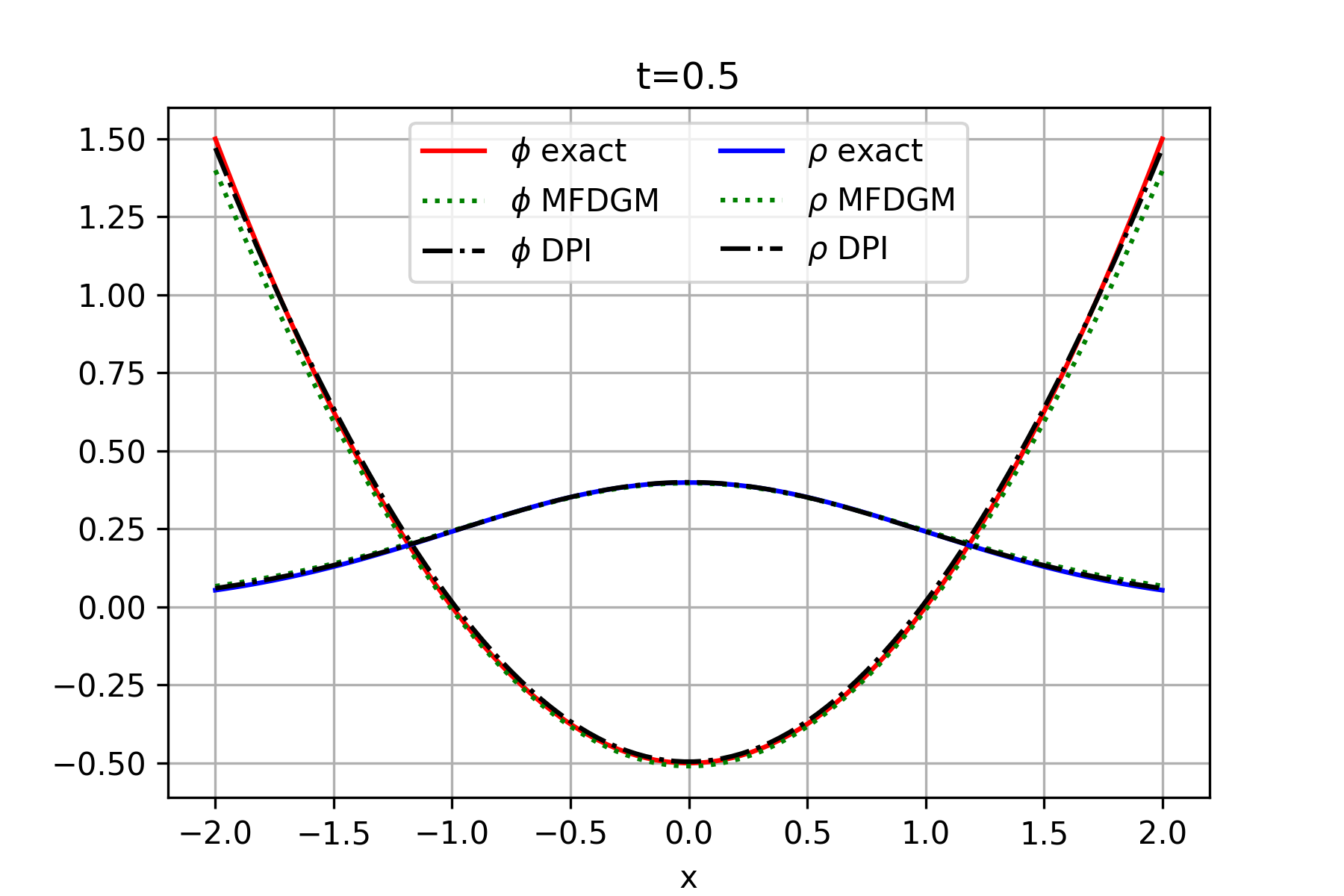}
\centering\includegraphics[width=3.8cm]{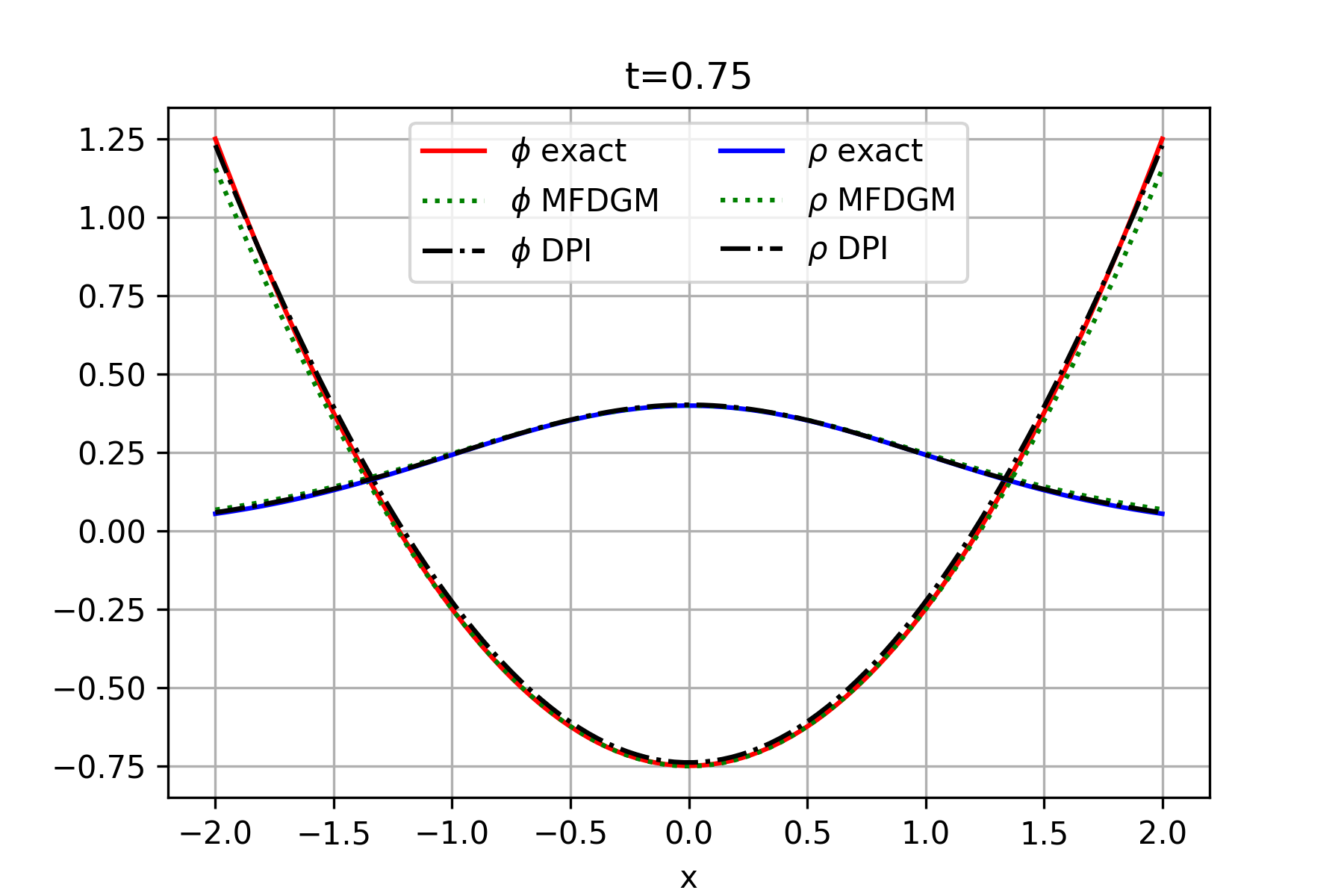}
\caption{We compare the exact solution of (\ref{test}) with the predicted solutions by DPI and MFDGM in one dimension at time points $( t = 0.25, 0.5, 0.75 )$ without congestion (\(\gamma = 0\)). The figures demonstrate that DPI and MFDGM predictions closely align with the exact solutions.}

\label{fig1}
\end{figure}

\begin{figure}[htbp]
\centering\includegraphics[width=5.8cm]{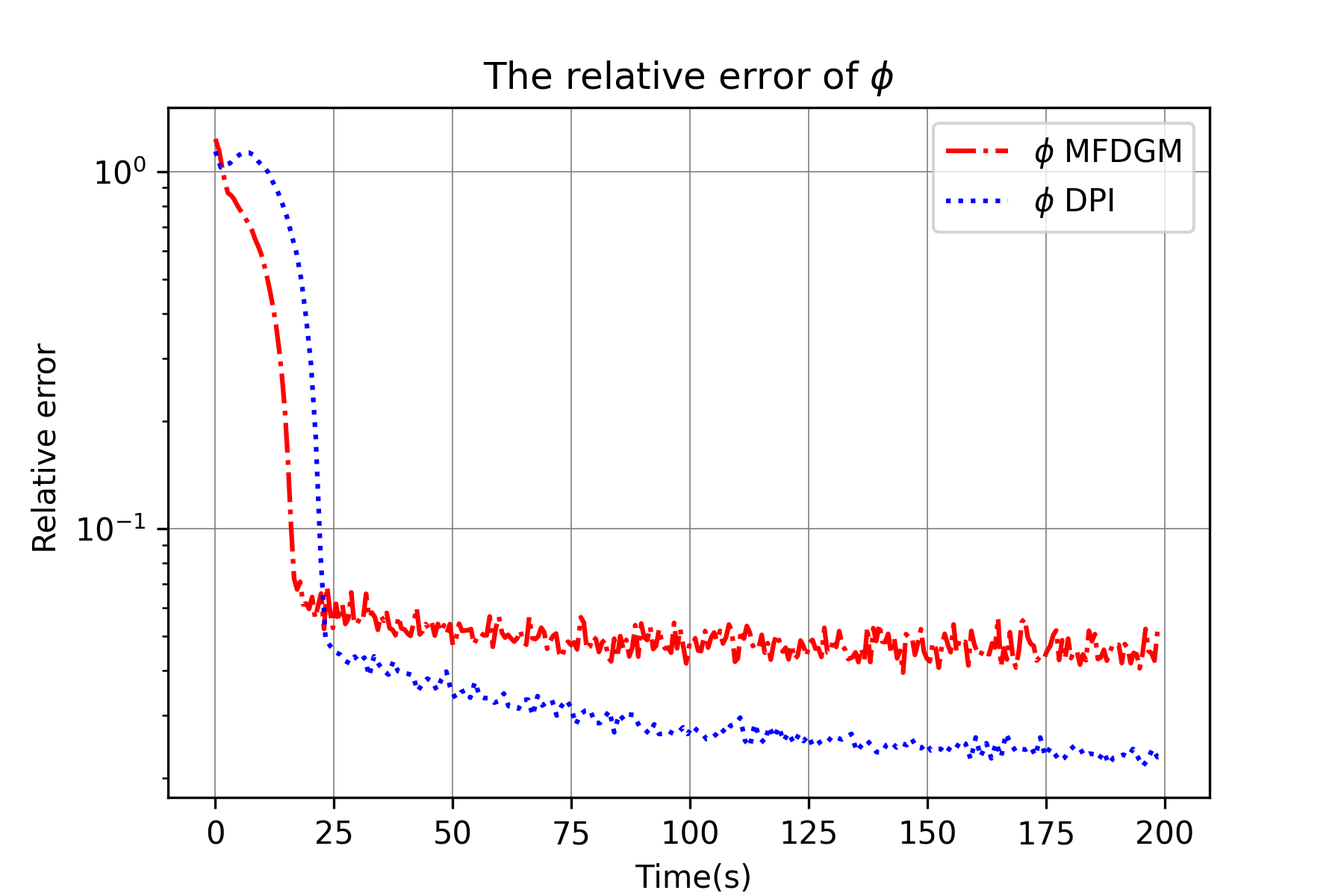}
\centering\includegraphics[width=5.8cm]{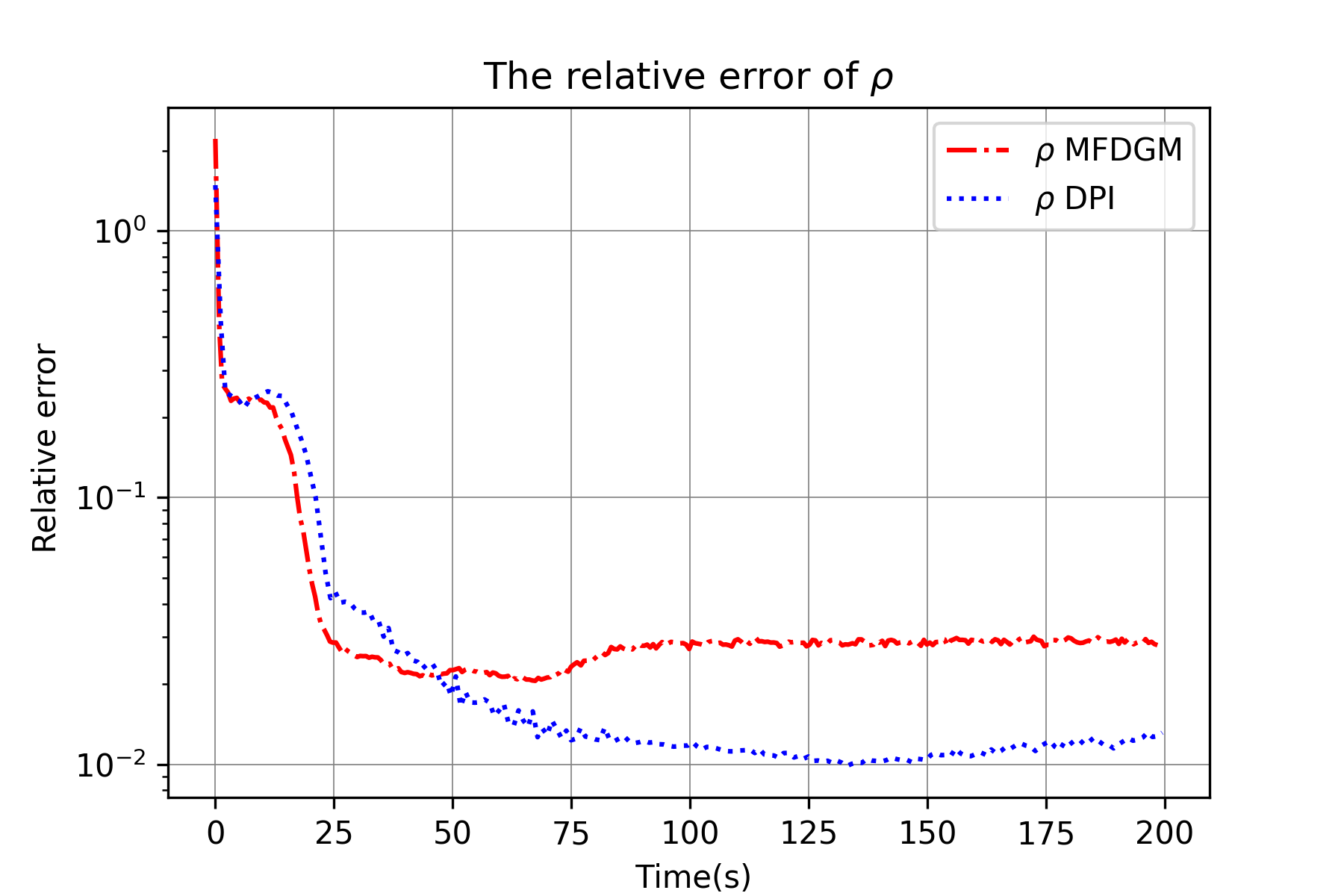}
\caption{We compare the relative error for $\phi$ and $\rho$ using two methods with $\gamma=0$. The figure demonstrates that DPI performs slightly better than MFDGM.} 
\label{fig2}
\end{figure}

\begin{figure}[htbp]
\centering\includegraphics[width=5.8cm]{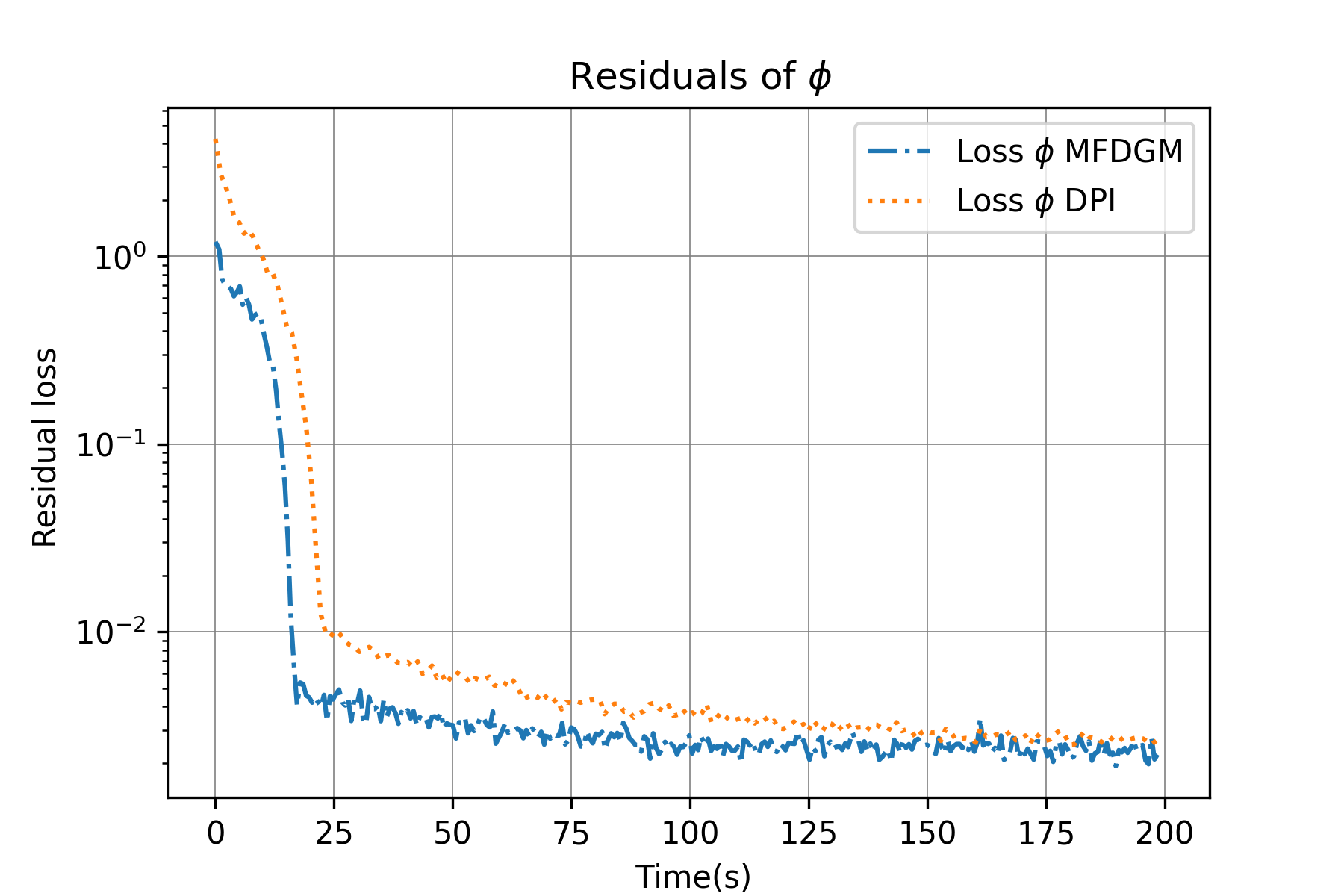}
\centering\includegraphics[width=5.8cm]{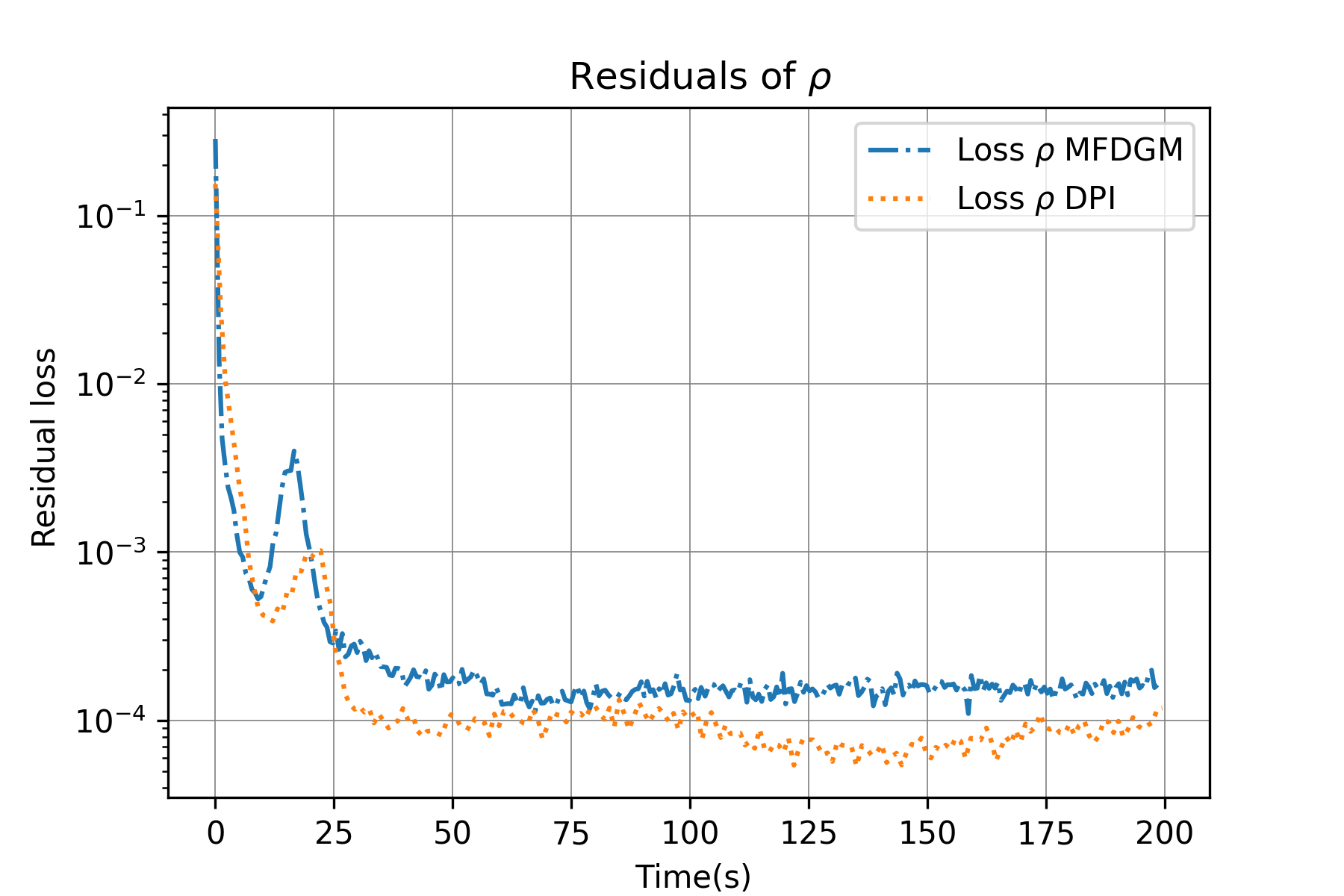}
\caption{To monitor the convergence, we show the residual losses for $\phi$ and $\rho$ using two methods with $\gamma=0$.}
\label{fig3}
\end{figure}

{\bf Test 1:} Assuming a congestion-free scenario ($\gamma=0$) in a one-dimensional system of partial differential equations \ref{test}, we utilize Algorithm \ref{alg1} to obtain results by using a minibatch size of 50 samples at each iteration. Our approach employs neural networks with one hidden layer consisting of 100 neurons each, with Softplus activation function for $N_{\omega}$ and Tanh activation function for  $N_{\theta}$ and  $N_{\tau}$. To train the neural networks, we use ADAM with a learning rate of $10^{-4}.$ and weight decay of $10^{-3}.$ We adopt ResNet as the architecture of the neural networks, with a skip connection weight of 0.5. We kept the same parameters for the second method MFDGM. The numerical results are presented in Figure \ref{fig1}, where we compare the approximate solutions obtained by two methods to the exact solutions at different times.
To assess the performance of the methods, we compute the relative error between the model predictions and the exact solutions on a $100 \times 100$ grid within the domain $[0,1]\times[-2,2]$, see Figure \ref{fig2}. Furthermore, we monitored the convergence of our approach by plotting the two residual losses, as defined in Algorithm \ref{alg1} and the MFDGM algorithm, in Figure \ref{fig3}. Optimal selection of architecture and hyperparameters for the neural networks plays a crucial role in attaining the desired outcomes for both methods see experiments in \cite{ASSOULI2023113802}.
\\

\begin{figure}
\centering\includegraphics[width=3.8cm]{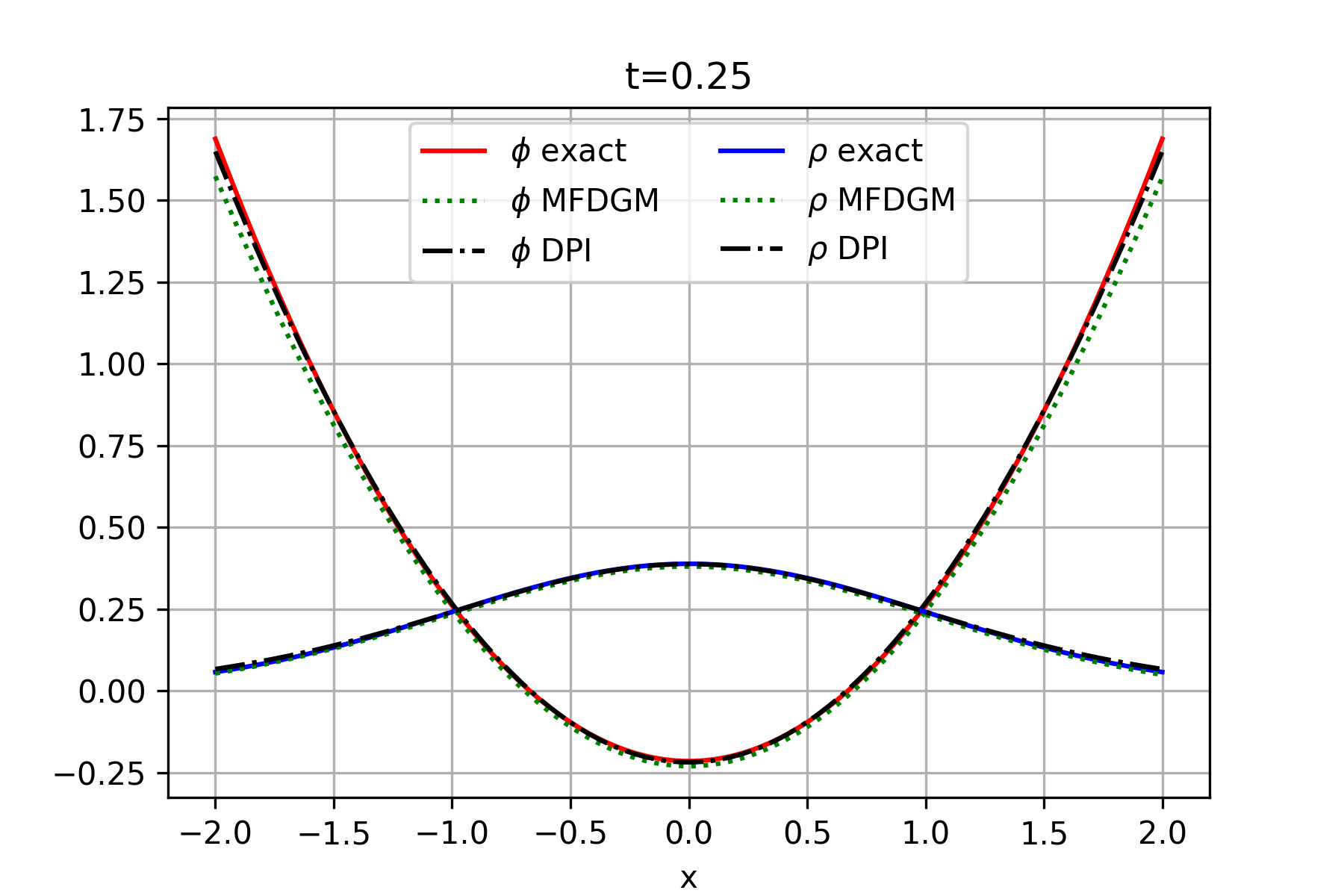}
\centering\includegraphics[width=3.8cm]{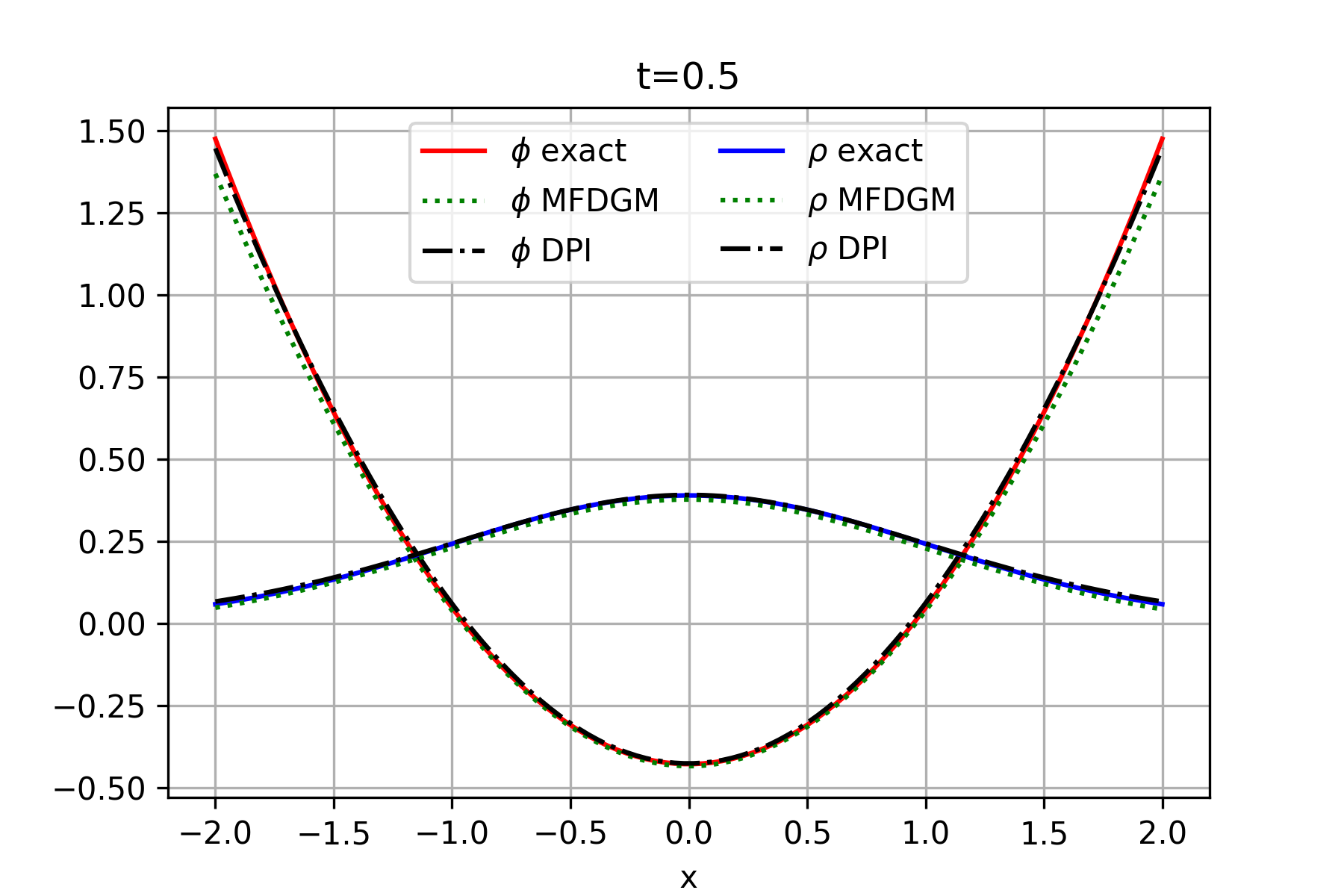}
\centering\includegraphics[width=3.8cm]{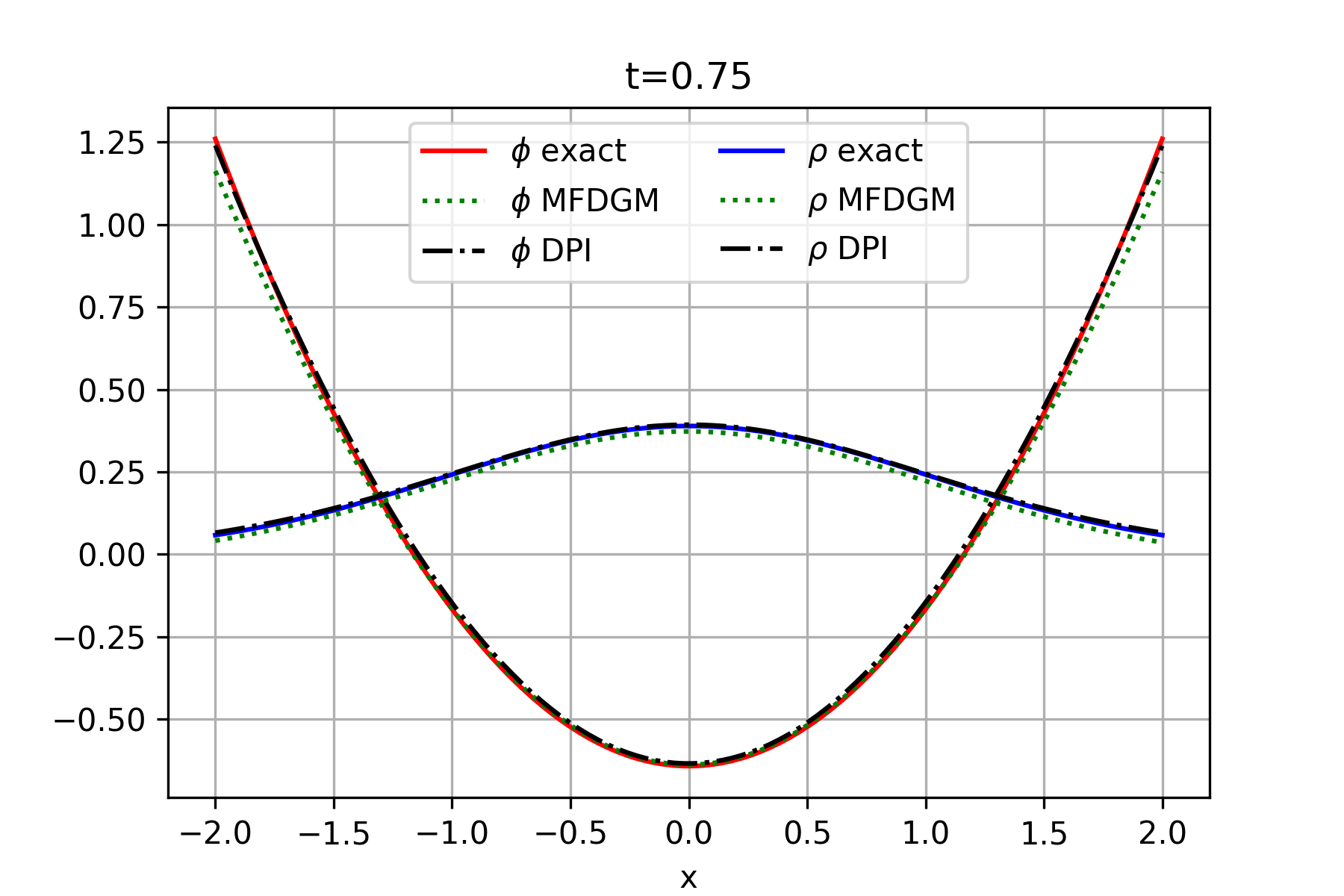}
\caption{We compare the exact solution of (\ref{test}) with the predicted solutions by DPI and MFDGM in one dimension at time points $( t = 0.25, 0.5, 0.75 )$, but now with congestion (\(\gamma = 0.1\)). The figures demonstrate that DPI and MFDGM predictions closely align with the exact solutions.} 
\label{fig4}
\end{figure}

\begin{figure}[htbp]
\centering\includegraphics[width=5.8cm]{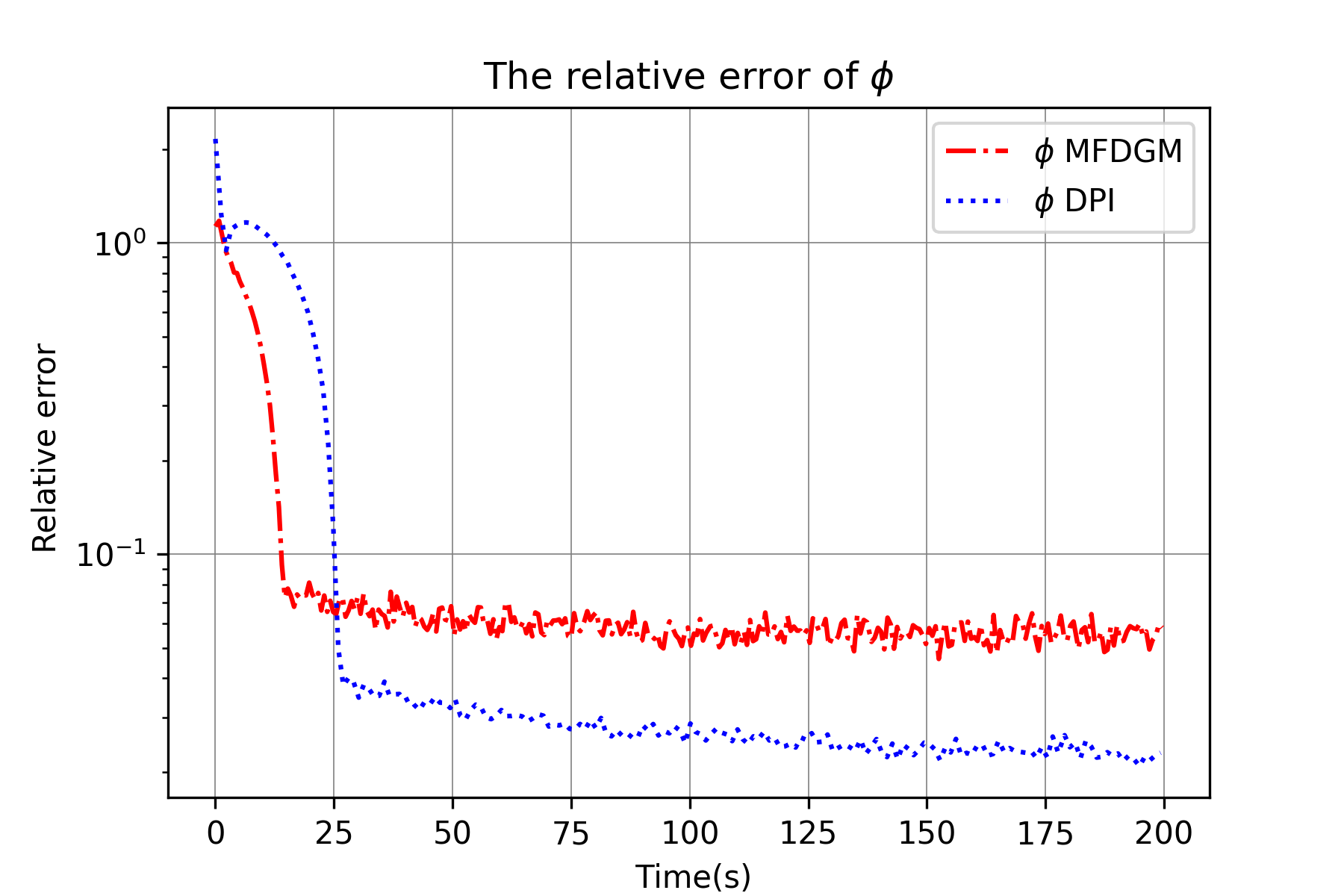}
\centering\includegraphics[width=5.8cm]{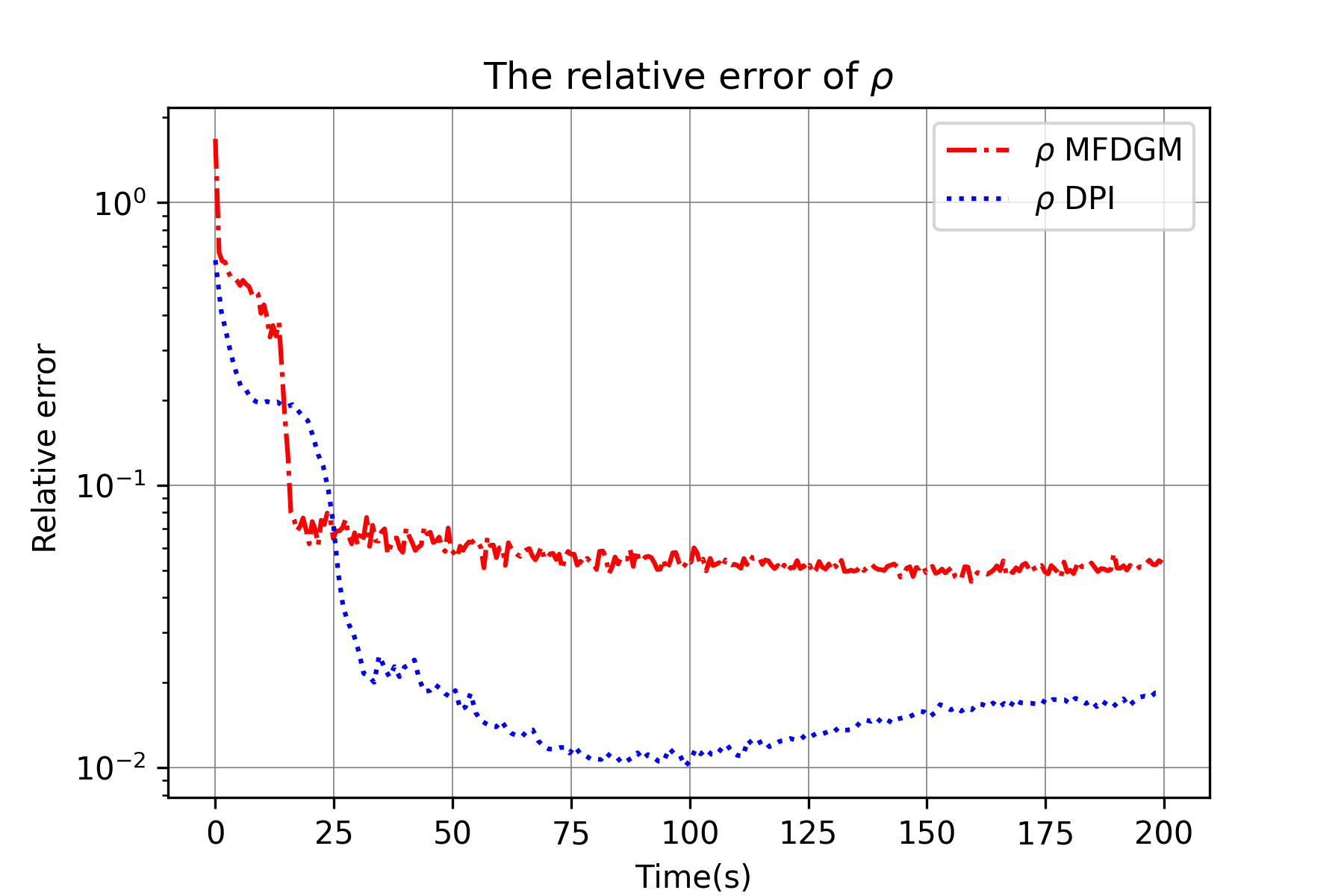}
\caption{Comparison of Relative Error for $\phi$ and $\rho$ using Two Methods with $\gamma=0.1$. The figure demonstrates again that DPI performs slightly better than MFDGM.} 
\label{fig5}
\end{figure}

\begin{figure}[htbp]
\centering\includegraphics[width=5.8cm]{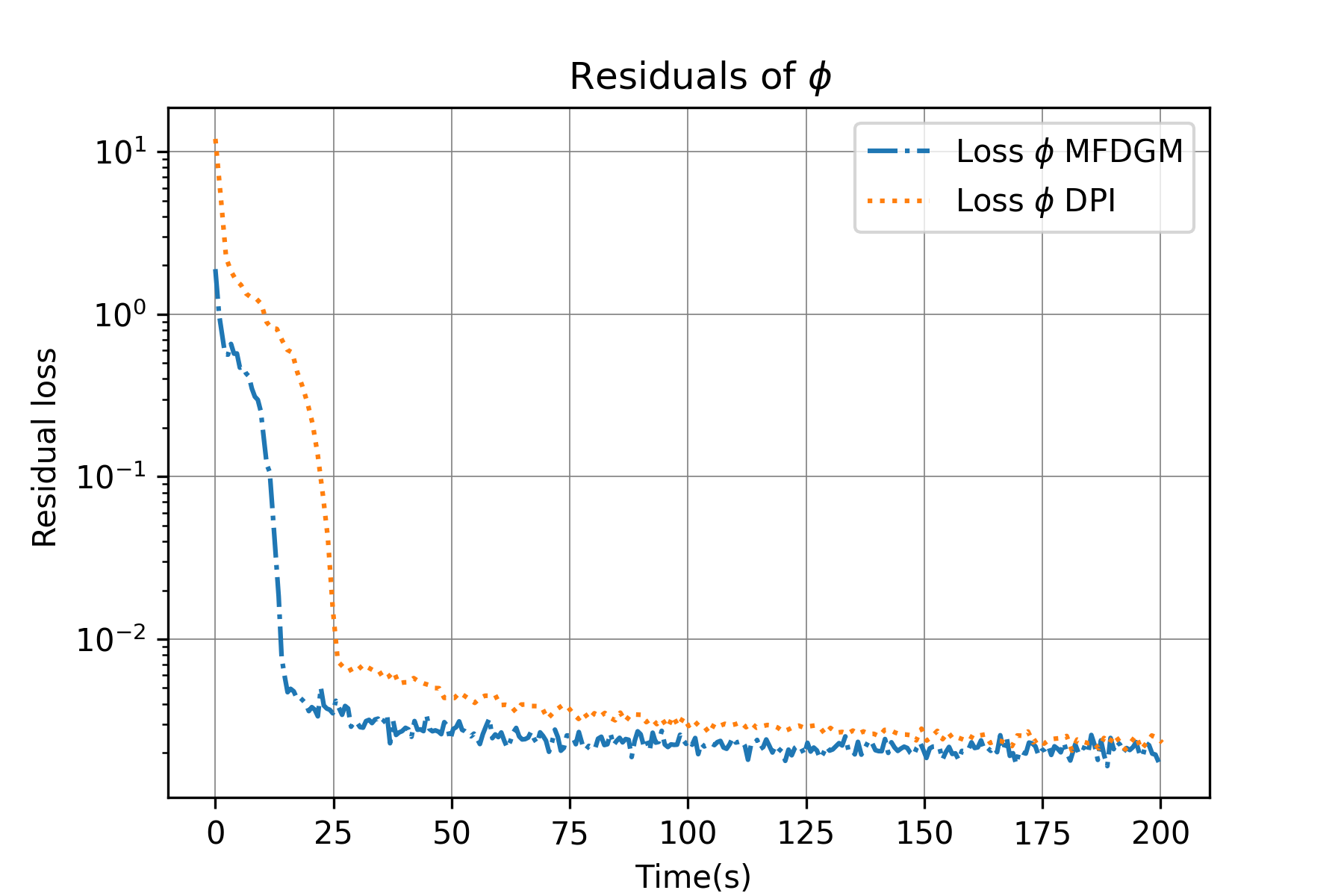}
\centering\includegraphics[width=5.8cm]{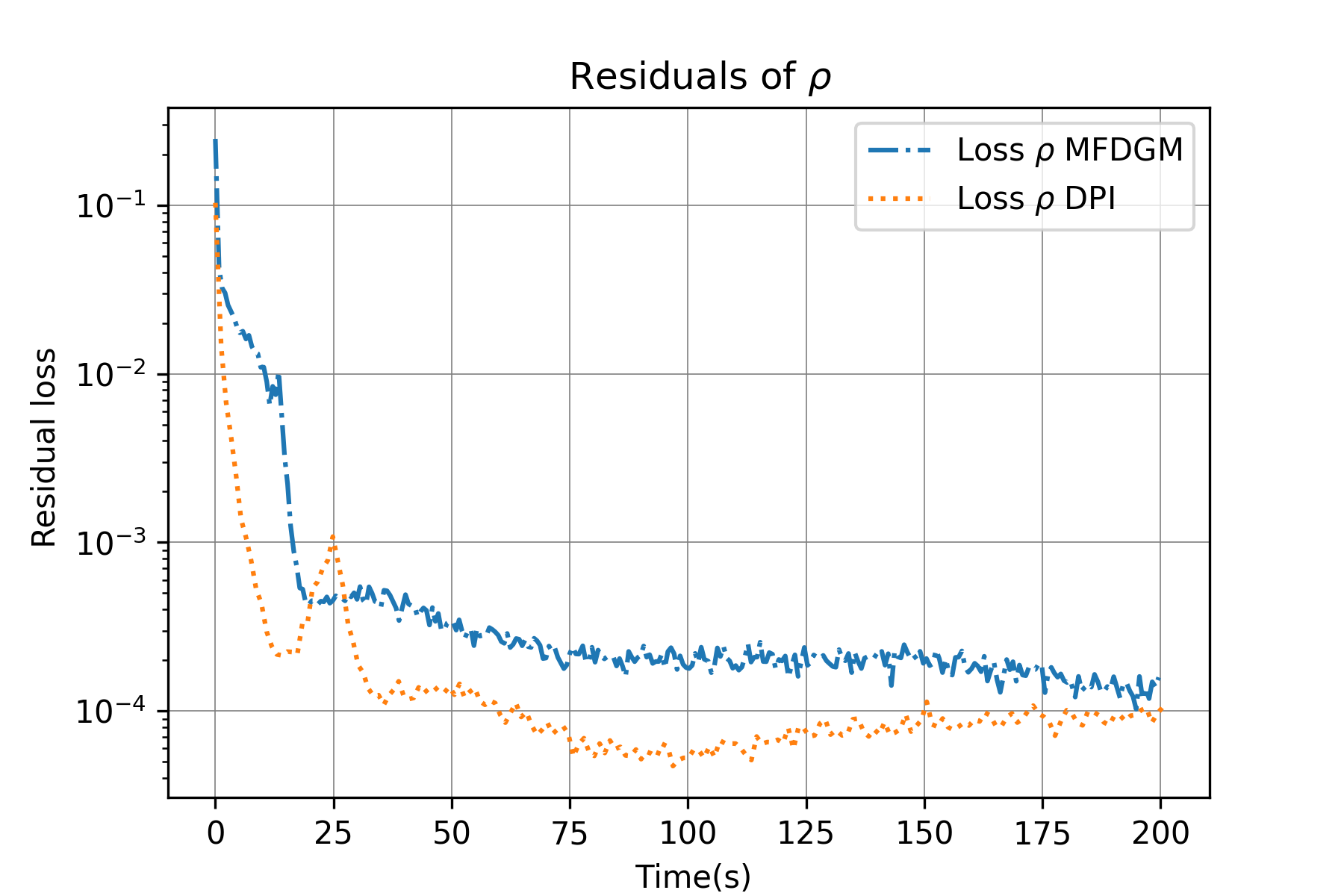}
\caption{To monitor the convergence, we show the residual losses for $\phi$ and $\rho$ using two methods with $\gamma=0.1$.}
\label{fig6}
\end{figure}

{\bf Test 2:} In order to investigate the impact of congestion on the previous case, we repeat the same experiment with the same parameters, but this time we consider a non-zero congestion parameter $(\gamma=0.1)$. We present the numerical results in Figures \ref{fig4}, \ref{fig5}, \ref{fig6}.\\

In comparison to the MFDGM algorithm, DPI demonstrates a similar level of convergence, particularly noticeable when applied to congestion scenarios. As seen in Figures 3 and 6, the residual losses for both methods ultimately converge to the same point. The faster convergence of MFDGM is attributed to the use of only two neural networks trained simultaneously, as opposed to the three used in DPI.
However, it is important to note that loss alone is not sufficient to determine the convergence of a method. As illustrated in Figures 2 and 5, the relative error of DPI is actually better than that of MFDGM. This indicates that DPI provides more accurate results despite the slower convergence in terms of residual loss. This convergence advantage arises from DPI’s capacity to reduce the complexity of Partial Differential Equations (PDEs) through the integration of the PI algorithm. Furthermore, to address errors similar to those of MFDGM, we enforce constraints on DPI’s training process. This additional measure ensures enhanced accuracy in convergence. In contrast to MFDGM, which relies on monitoring convergence only through residual loss, DPI controls convergence based on the fixed-point method. This offers a more direct and robust means of ensuring convergence, especially where the analytic solution is unknown, as in the upcoming example.
\\

{\bf Test 3:} In this test, we aim to solve (\ref{test}) using PI. To enable the application of this method, we incorporate periodic boundary conditions into (6).  In their research, the authors utilized Policy Iteration methods to solve a Partial Differential Equation (PDE) system through finite-difference approximation. They adopted a uniform grid $\mathcal{G}$ and centred second-order finite differences for the discrete Laplacian while computing the Hamiltonian and divergence term in the FP equation using the Engquist-Osher numerical flux for conservation laws. The symbol $ \sharp$ represented the linear differential operators at the grid nodes. Specifically, in the context of one dimension, they implemented a uniform discretization of $E$ with $I$ nodes $x_i = ih$, where $h = 1/I$ is the space step. To discretize time, they employed an implicit Euler method for both the time-forward FP equation and the time-backward HJB equation, with a uniform grid on the interval $[0, T]$ with $N + 1$ nodes $t_n = n t$, where $t = T/N$ was the time step. To prevent any confusion, we adhere to the previously established notation. Therefore, we use $\phi$, $\rho$, and the policy $q$ to represent the vectors that approximate the solution on $\mathcal{G}$, and by $\phi_n, \rho_n,$ and $q_n$ the vectors on $\mathcal{G}$ approximating the solution and policy at time $t_n$. The algorithm we will use for the fully discretized System \ref{test} is as follows: we start with an initial guess  $q^{(0)}_n: \mathcal{G} \to \mathbb{R}^{2d}$ for $n = 0, . . . , N,$ and initial and final data $\rho_0, \phi_N: \mathcal{G} \to \mathbb{R}$. We then iterate on $k \geq 0$ until convergence,\\
(i) Solve for $n=0, \ldots, N-1$ on $\mathcal{G}$
$$
\left\{\begin{array}{l}
\rho_{n+1}^{(k)}-d t\left( \Delta_{\sharp} \rho_{n+1}^{(k)}+\operatorname{div}_{\sharp}\left(\rho_{n+1}^{(k)} q_{n+1}^{(k)}\right)\right)=\rho_n^{(k)} \\
\rho_0^{(k)}=\rho_0
\end{array}\right.
$$
(ii) Solve for $n=N-1, \ldots, 0$ on $\mathcal{G}$
$$
\left\{\begin{aligned}
\phi_n^{(k)} & -d t\left( \Delta_{\sharp} \phi_n^{(k)}-q_{n, \pm}^{(k)} \cdot D_{\sharp} \phi_n^{(k)}\right) \\
& =\phi_{n+1}+d t\left(\frac{1}{2}\left|q_{n+1, \pm}^{(k)}\right|^2+ x_i^2 +\gamma Ln\left(\rho_{n+1}^{(k)}\right)\right) \\
\phi_N^{(k)} & =\phi_N
\end{aligned}\right.
$$
(iii) Update the policy $q_n^{(k+1)}=D_{\sharp} \phi_n^{(k)}$ on $\mathcal{G}$ for $n=0, \ldots, N$, and set $k \leftarrow$ $k+1$.
In the following test, we choose  $\gamma=0$,  T = 1 for the final time, and $K=50$. The grid consisted of $I=200$ nodes in space and $N=200$ nodes in time. The initial policy was initialized as $q^{(0)}_n \equiv (0, 0)$ on $\mathcal{G}$ for all $n$.
\begin{figure}
\centering\includegraphics[width=3.8cm]{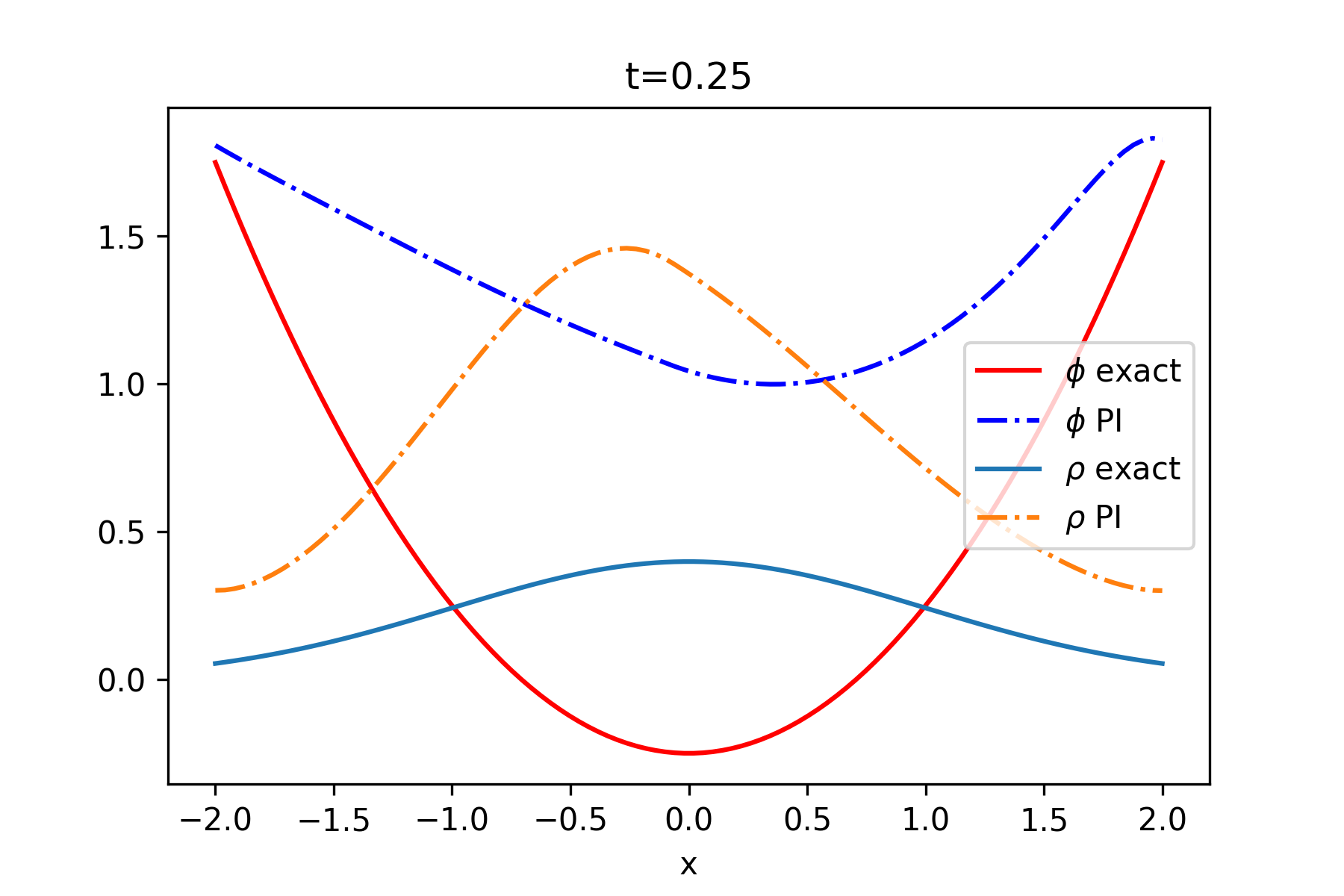}
\centering\includegraphics[width=3.8cm]{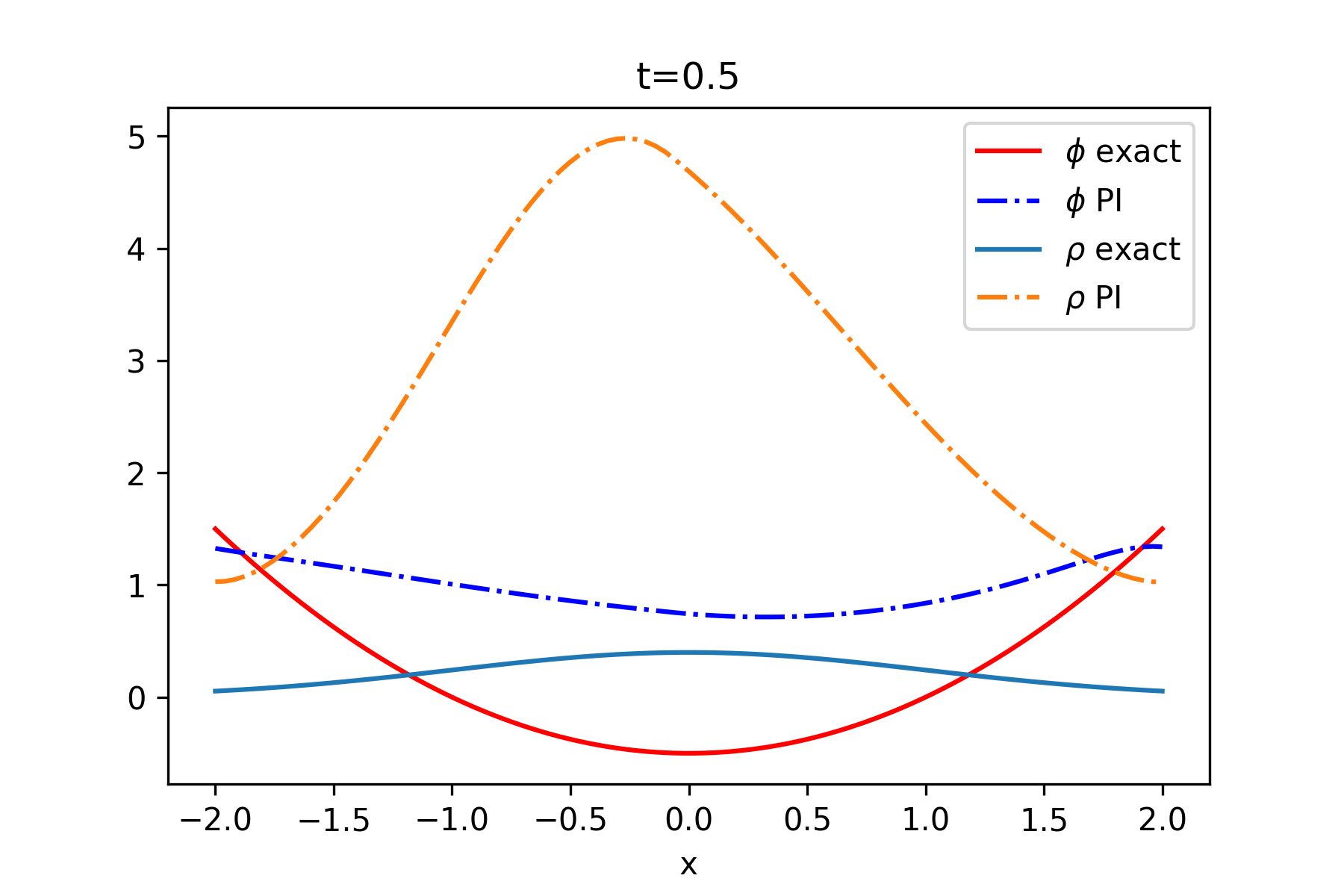}
\centering\includegraphics[width=3.8cm]{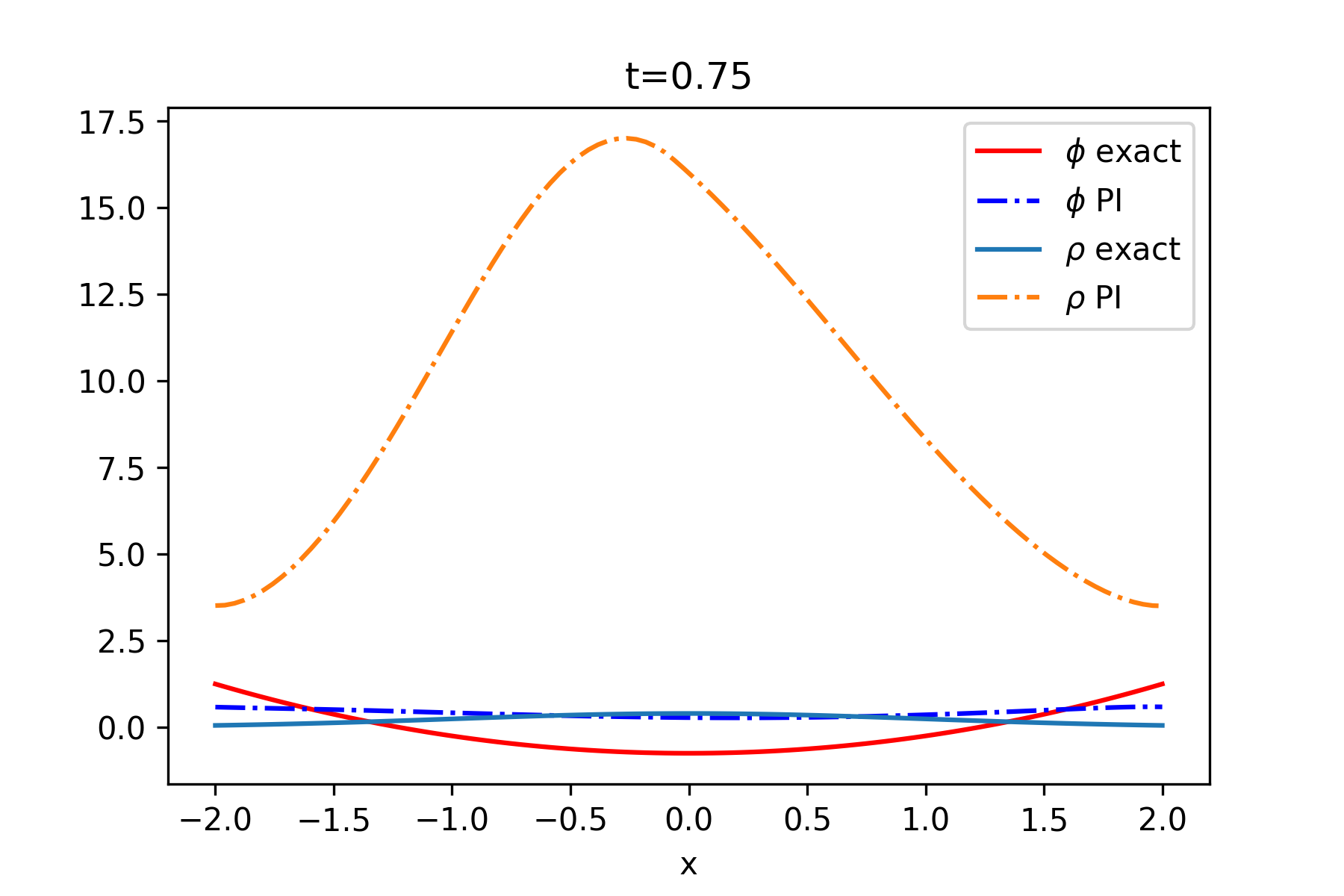}
\caption{We compare the exact solution (1) with the predicted solution by PI  in one dimension at different time points \( t = 0.25, 0.5, 0.75 \) under no congestion (\(\gamma = 0\)). The figures show that PI  prediction does not closely match the exact solutions, suggesting that PI may not be effective for separable cases.}
\label{fig111}
\end{figure}\\

The results of the study Figure \ref{fig111} showed that policy iteration may not be effective in solving MFG problems in separable cases. However, the use of deep learning instead of finite difference methods has exhibited promise in extending the applicability of this approach to MFG in the separable case. \\

\begin{figure}
\centering\includegraphics[width=5.8cm]{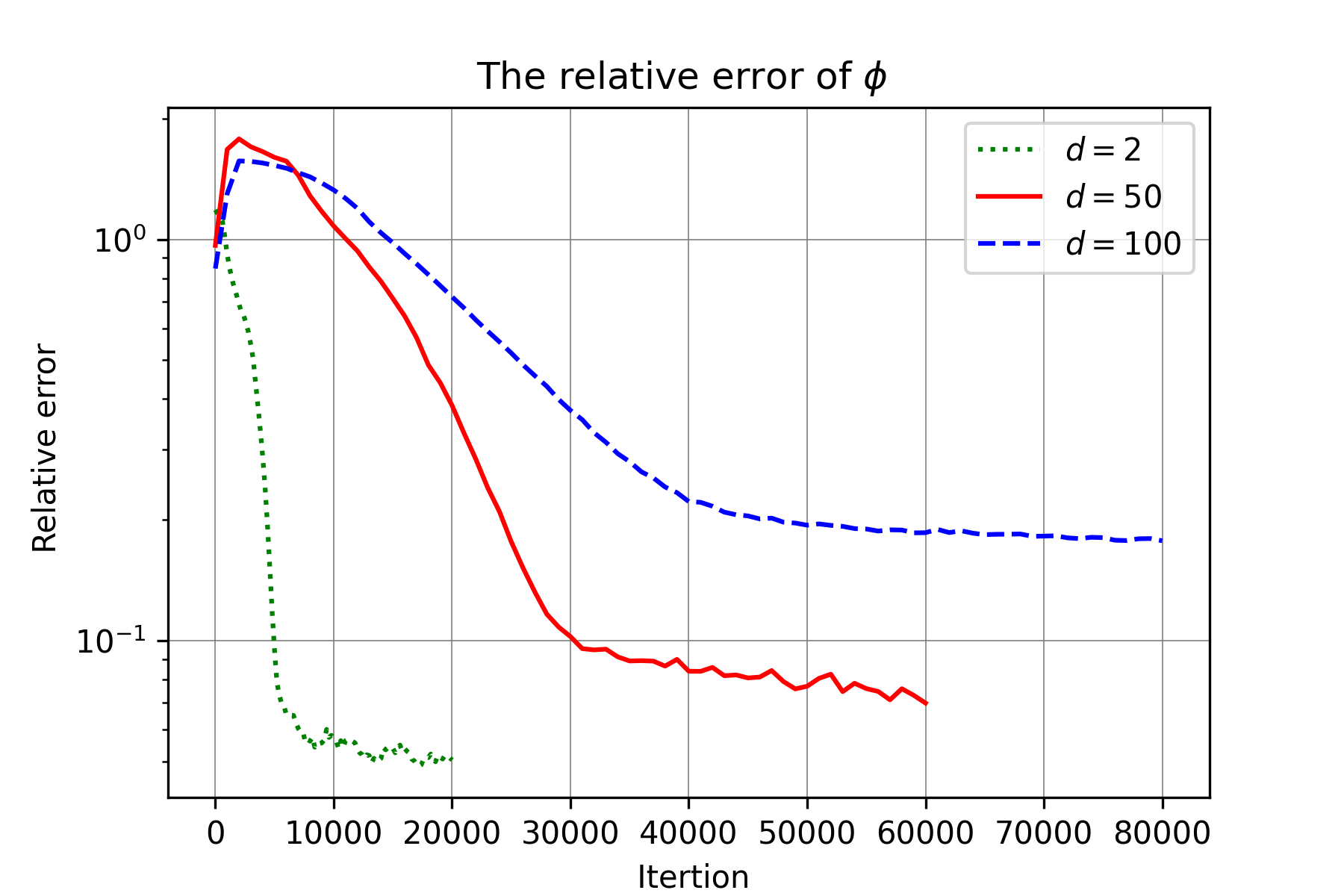}
\centering\includegraphics[width=5.8cm]{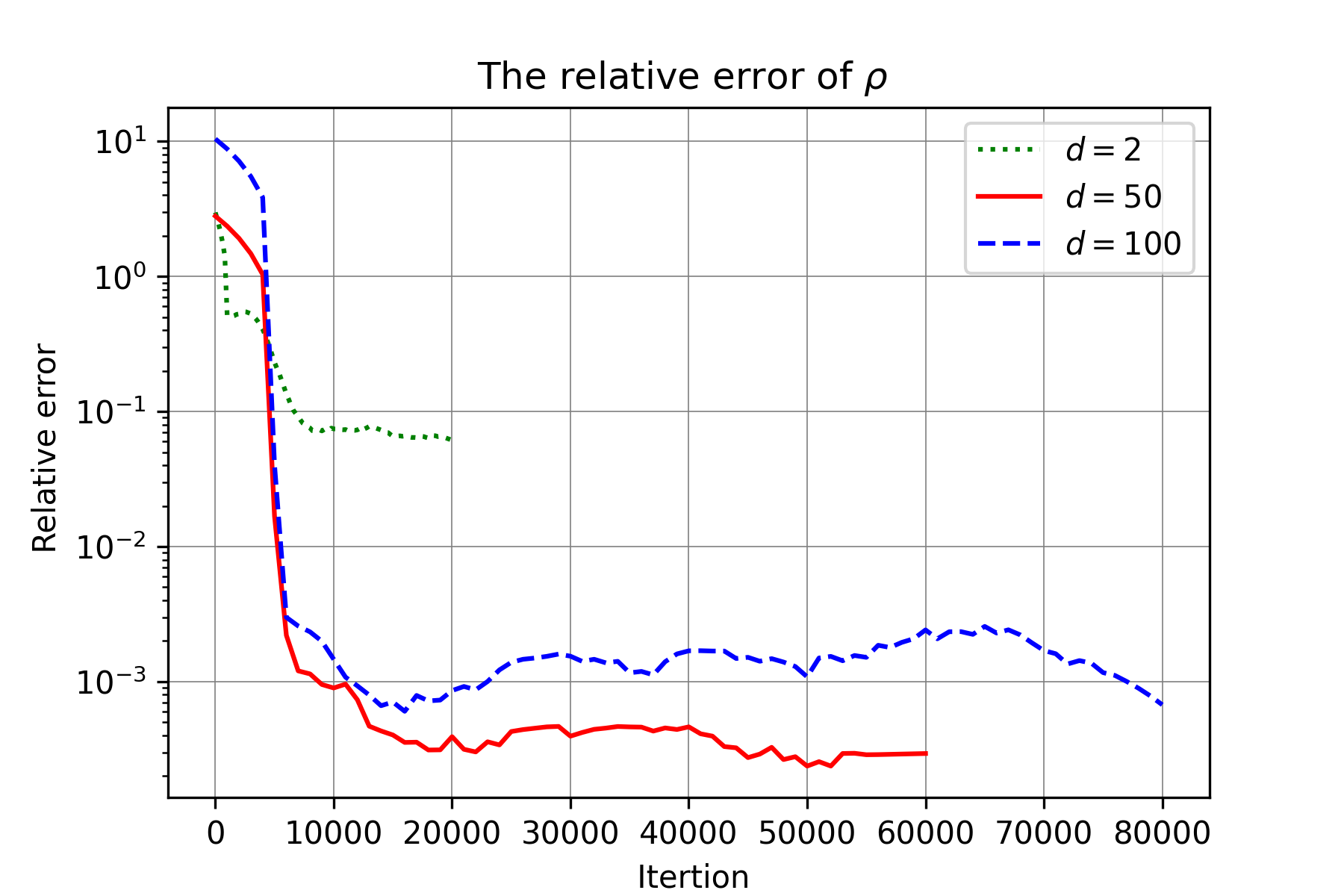}
\caption{The relative error for $\phi$ and $\rho$ of the DPI method is examined across high-dimensional cases, ranging from dimensions 
$d=2$ to $d=50$ and 
$d=100$. The figure illustrates that DPI can approximates solutions even in the challenging scenario of 100 dimensions.}
\label{100}
\end{figure}

{\bf Test 4:} The effectiveness of our approach in high-dimensional cases is tested in this experiment and the next section. We solve the MFG system (\ref{test}) for dimensions 2, 50, and 100 and present the results in Figure \ref{100}. To train the neural networks, we use a minibatch of 100, 500, and 1000 samples respectively for d=2, 50 and 100. The neural networks have a single hidden layer consisting of 100, 200 and 256 neurons each respectively. We utilize the Softplus activation function for $N_{\omega}$ and the Tanh activation function for $N_{\tau}$ and $N_{\theta}$. We use ADAM with a learning rate of $10^{-4}$, weight decay of $10^{-4}$, and employ ResNet as the architecture with a skip connection weight of 0.5. 

The results in Figure \ref{100} were obtained after applying a Savgol filter to enhance the clarity of the curves. This result highlights how DPI can handle many dimensions by using the neural network's strengths, which in turn improves the usefulness of the PI algorithm. We will explore more complex examples in the next section.

\subsection{High-Dimensional}
In this section, our goal is to further assess the effectiveness of our approach in high-dimensional contexts by examining two specific examples characterized by non-separable Hamiltonians, as previously introduced in \cite{lauriere2021policy}.\\

\begin{figure}
\centering\includegraphics[width=3.8cm]{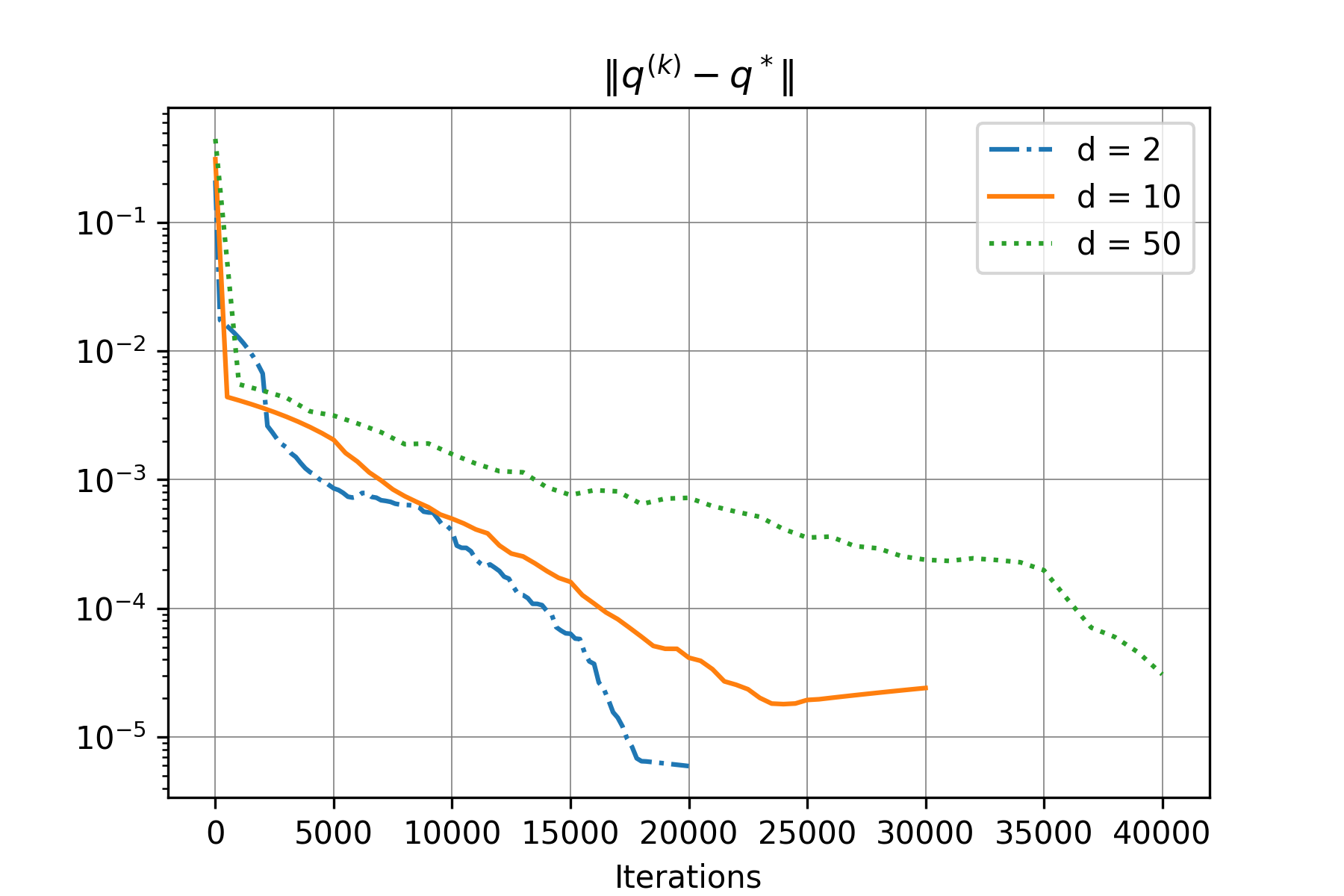}
\centering\includegraphics[width=3.8cm]{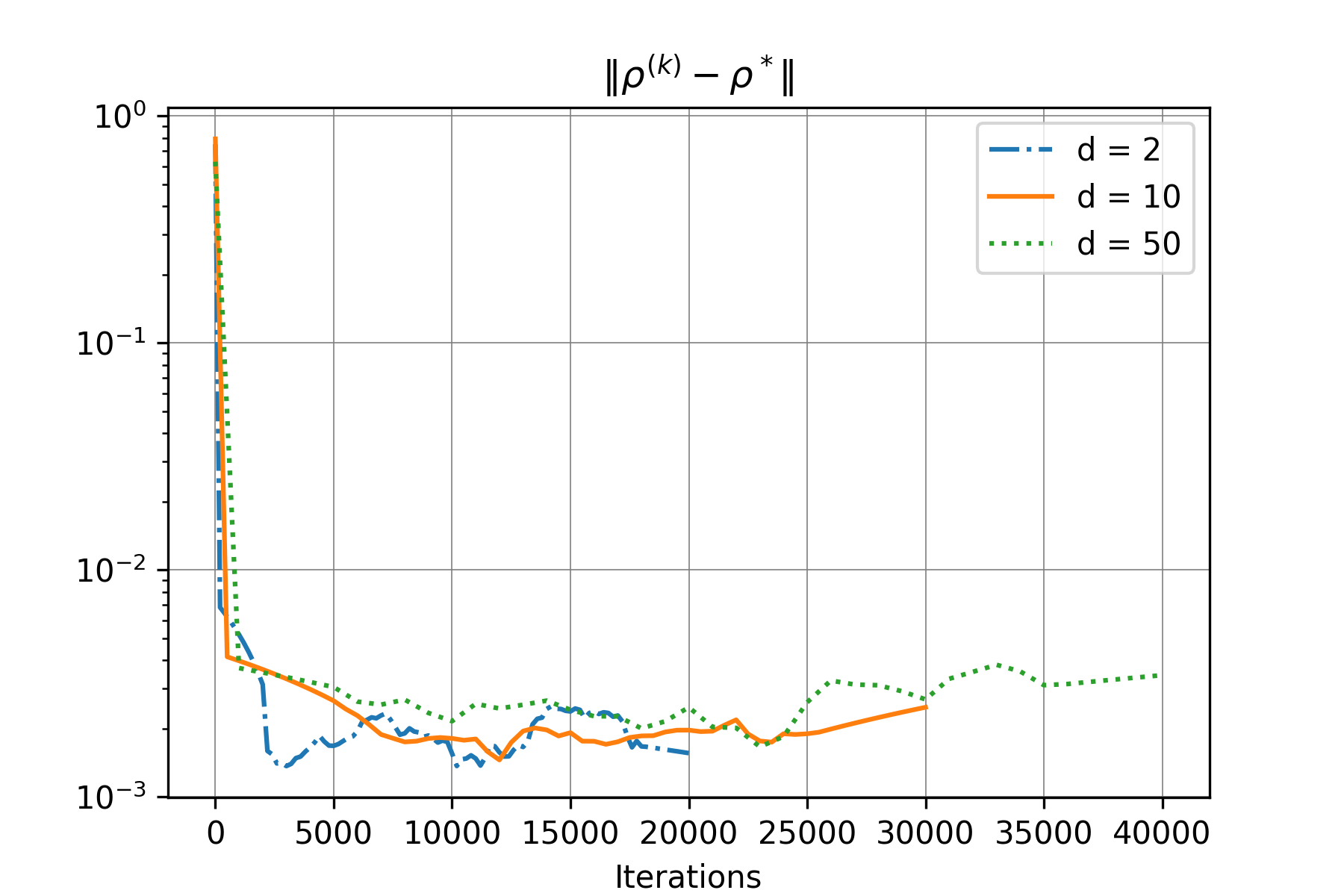}
\centering\includegraphics[width=3.8cm]{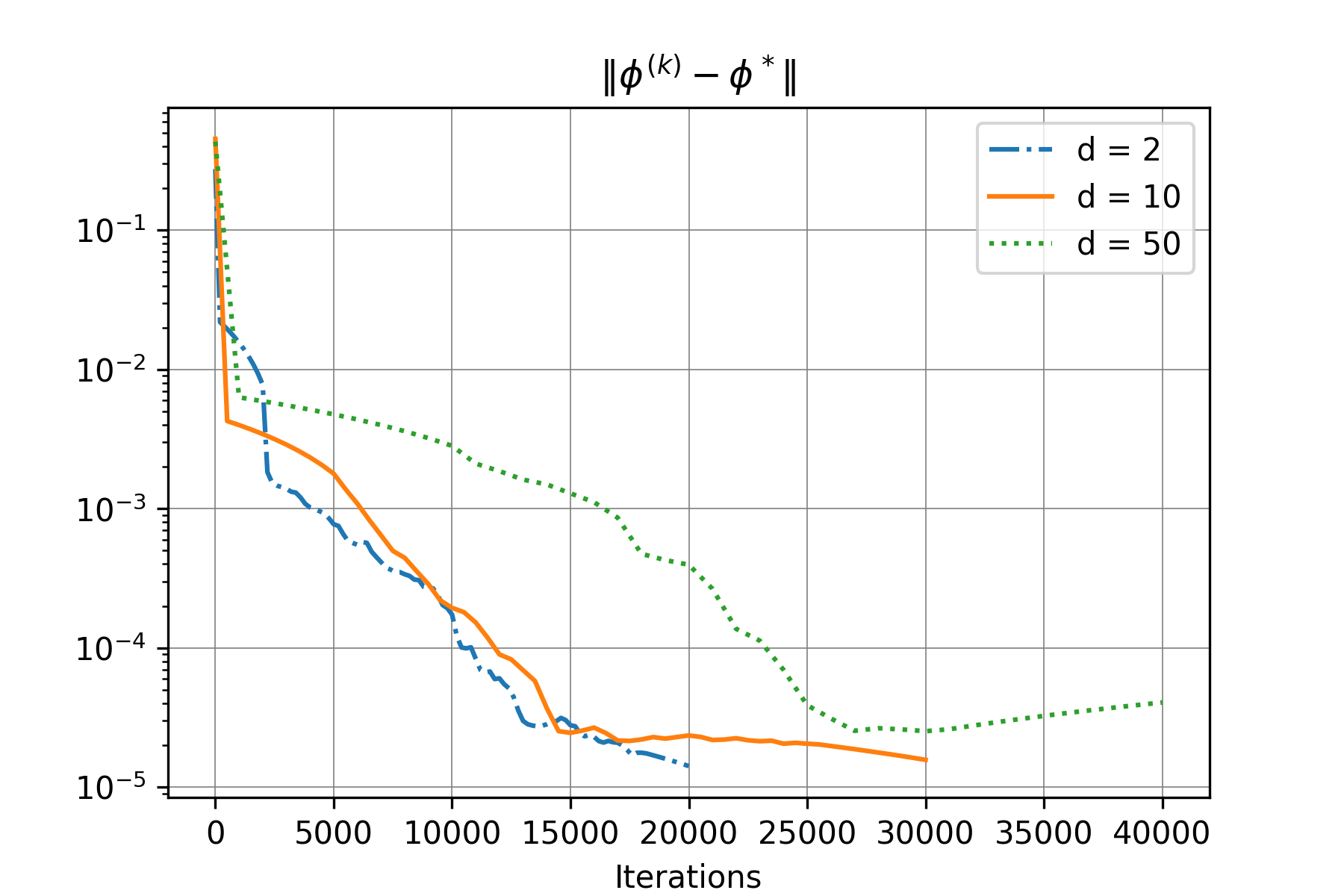}
\caption{Example1: The $L^{\infty}$ Distance between $\rho^{(k)}$,  $\phi^{(k)}$  and $q^{(k)}$ from DPI and the final solution $\rho^{*}$,  $\phi^{*}$  and $q^{*}$ from fixed point algorithm for  d=2, 10 and 50. }
\label{figex1}
\end{figure}

{\bf Example 1:} We consider the following problem in the stochastic case with $\nu = 0.3$. The problem is defined within the domain $\Omega=[0,1]^d$, with a fixed final time of $T = 1$. In this context, the terminal cost is specified as $g=0,$ and the initial density is characterized by a Gaussian distribution.\\
The corresponding MFG system is,
\begin{equation}\label{exp1}
\left\{
\begin{array}{rrrrr}
-\phi_t-\nu \Delta \phi+\frac{1}{2(1+4\rho)}||\phi_x||^2=0,\\
 \rho_t -\nu \Delta \rho-\operatorname{div}(\rho\frac{\phi_x}{1+4\rho})=0,\\
\rho(0,x)=\left(\frac{1}{2\pi}\right)^{d/2}e^{-\frac{ ||x-0.25||^2}{2}}, \\ \phi(T, x)=0. 
\end{array}
\right.
\vspace{-0.25cm}
\end{equation}\\
We addressed the MFG system (\ref{exp1}) across three different dimensions: 2, 10, and 50. The results of this analysis are presented in Figure \ref{figex1}. For the training of neural networks, we adopted varying minibatch sizes: 100, 500, and 1000 samples, corresponding to $d= 2, 10$, and $100$. These neural networks consist of a single hidden layer with 100. We applied the Softplus activation function for $N_{\omega}$ and the Tanh activation function for $N_{\tau}$ and $N_{\theta}$. Our optimization approach used ADAM with a learning rate of $10^{-4}$ and a weight decay of $10^{-4}$. We employed the ResNet architecture with a skip connection weight of 0.5. To enhance the smoothness of the curves, the results underwent a Savgol filter.

\begin{figure}
\centering\includegraphics[width=3.8cm]{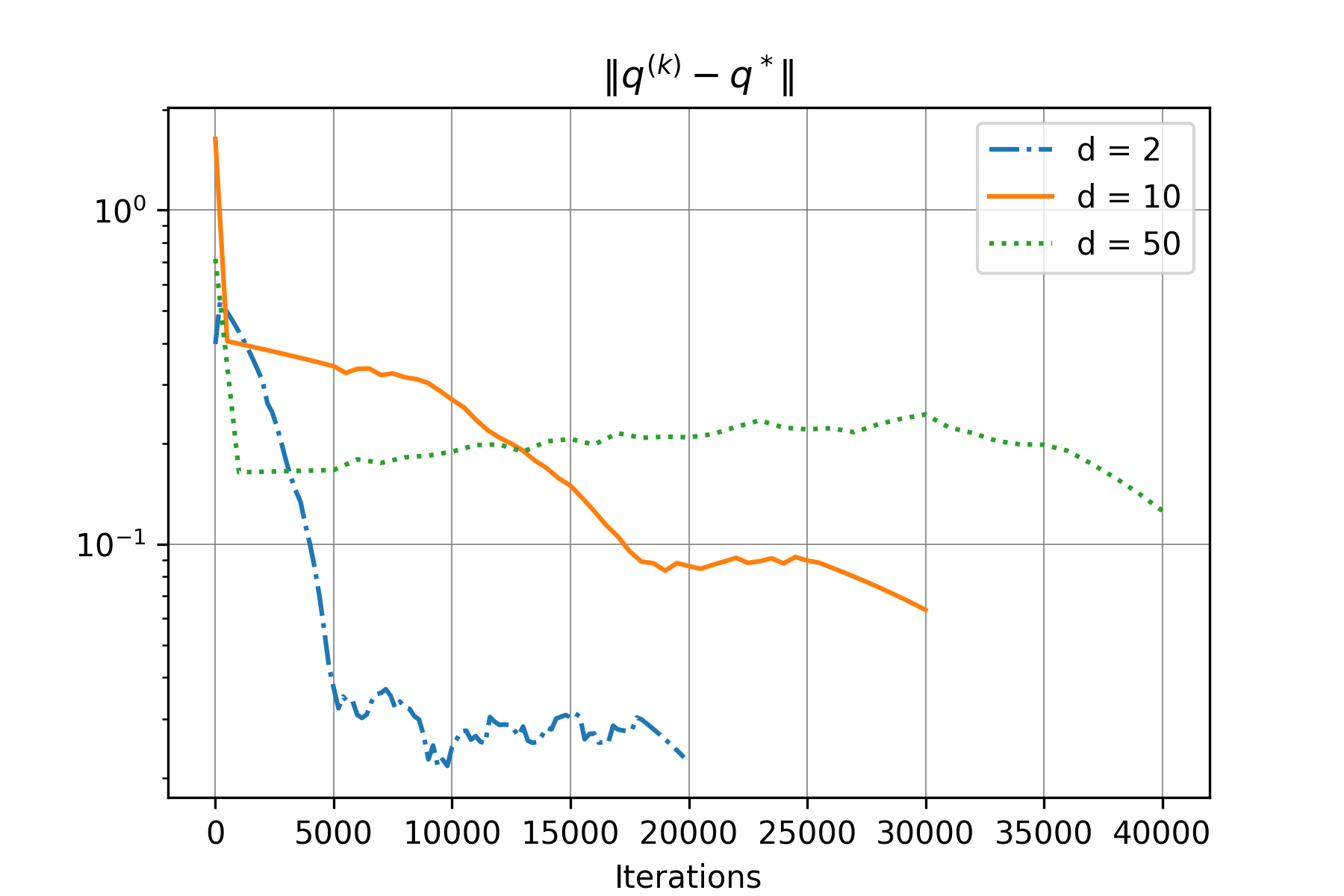}
\centering\includegraphics[width=3.8cm]{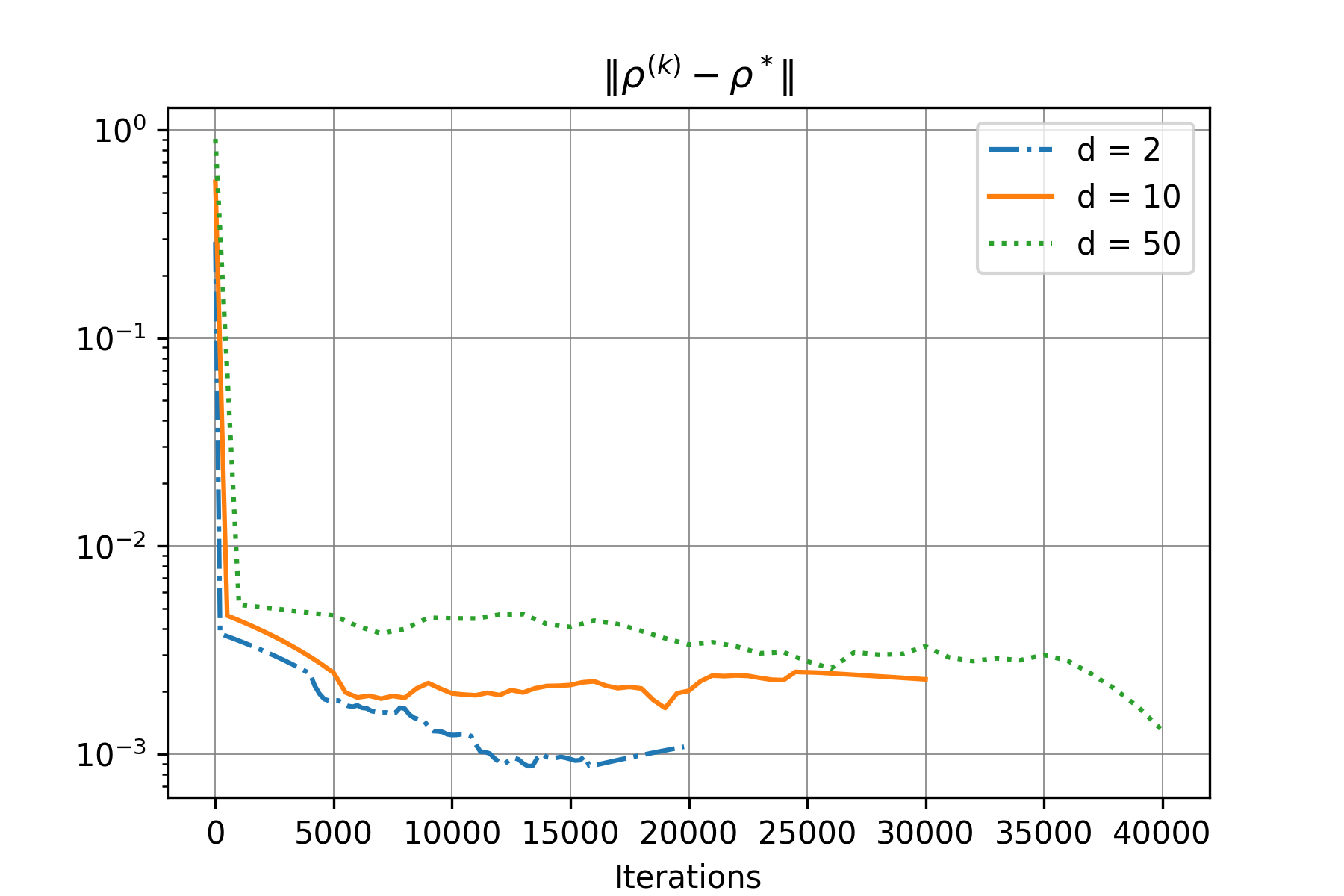}
\centering\includegraphics[width=3.8cm]{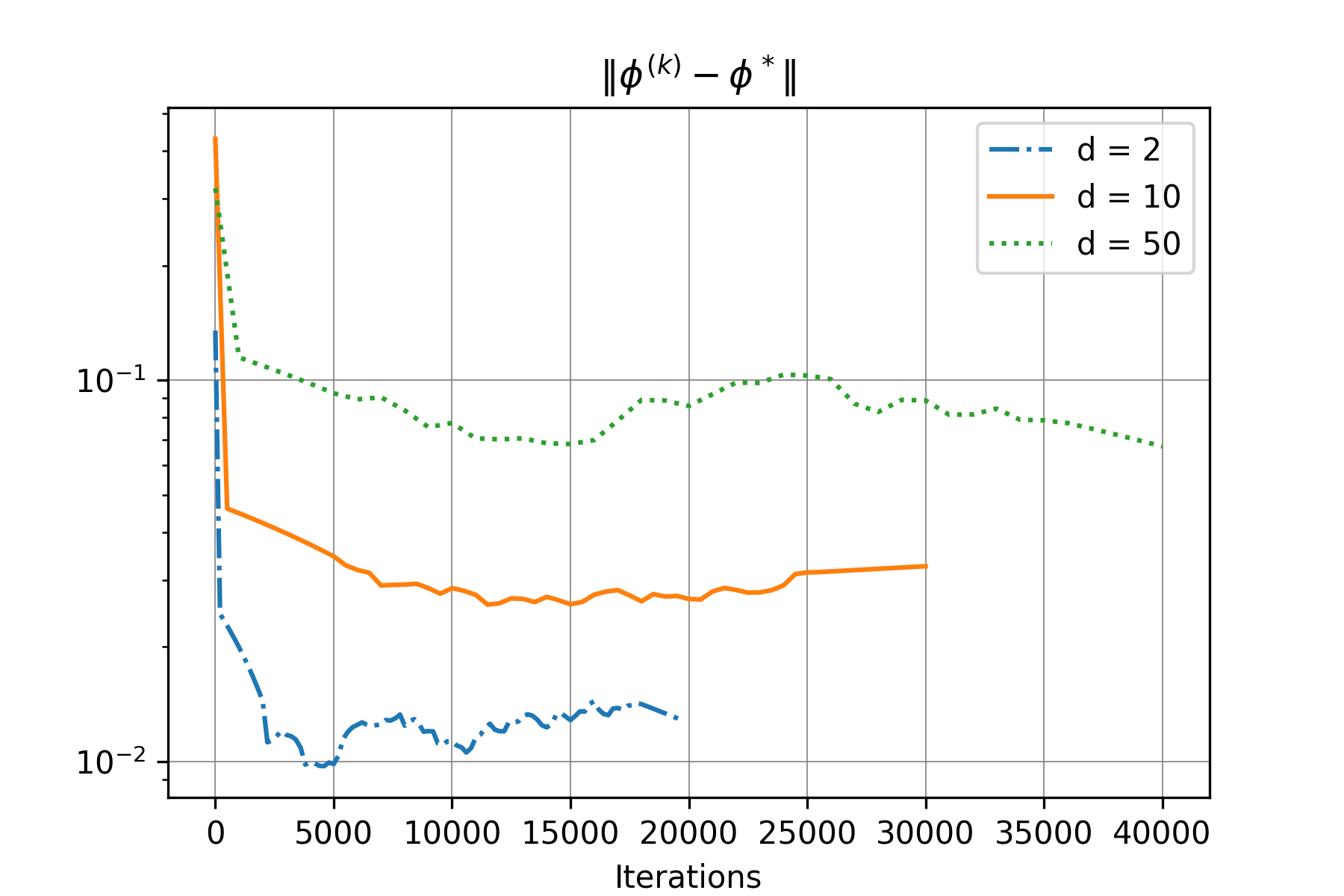}
\caption{Example2: The $L^{\infty}$ Distance between $\rho^{(k)}$,  $\phi^{(k)}$  and $q^{(k)}$ from DPI and the final solution $\rho^{*}$,  $\phi^{*}$  and $q^{*}$ from fixed point algorithm for  d=2, 10 and 50. } 
\label{figex2}
\end{figure}

{\bf Example 2:} Now we give the following problem with a terminal cost that motivates the agents to direct their movements toward particular sub-regions within the domain.\\
The corresponding MFG system is,
\begin{equation}\label{exp2}
\left\{
\begin{array}{rrrrr}
-\phi_t-\nu \Delta \phi+\frac{1}{2\rho^{1/2}}||\phi_x||^2=0,\\
 \rho_t -\nu \Delta \rho-\operatorname{div}(\rho\frac{\phi_x}{\rho^{1/2}})=0,\\
\rho(0,x)=\left(\frac{1}{2\pi}\right)^{d/2}e^{-\frac{ ||x-0.25||^2}{2}}, \\ \phi(T, x)=\sum_{i=1}^{d} cos(2\pi x_i). 
\end{array}
\right.
\vspace{-0.25cm}
\end{equation}\\
We perform the same process as previously described in Example 1 for the MFG system (\ref{exp2}), and the results of this analysis are presented in Figure \ref{figex2}.

In Figures \ref{figex1} and \ref{figex2}, we report results on the convergence with respect to the number of iterations. Specifically, we compare the \(L^{\infty}\) distance between \(\rho^{(k)}\), \(\phi^{(k)}\), and \(q^{(k)}\) computed by DPI, and the final solution \(\rho^{*}\), \(\phi^{*}\), and \(q^{*}\) from the fixed-point iteration algorithm, as described in \cite{lauriere2021policy}. This approach allows us to establish a benchmark solution when the exact solution is unknown, instead of monitoring convergence through residual losses as done in MFDGM.

\subsection{Trafic Flow}\label{tf}
To assess the DPI method, we conducted experiments on a traffic flow problem for an autonomous vehicle used to test the effectiveness of MFDGM. We specifically chose this problem because it is characterized by a non-separable Hamiltonian, making it a more complex and challenging problem. In our experiments, we applied the PI and DPI Methods to the problem and evaluated its performance in stochastic cases $\nu = 0.1$. We consider the traffic flow problem on the spatial domain $\Omega=[0,1]$ with dimension $d=1$ and final time $T=1$. The terminal cost $g$ is set to zero and the initial density $\rho_0$ is given by a Gaussian distribution, $\rho_0(x)=0.05-0.9 \exp\left(\frac{-1}{2}\left(\frac{x-0.5}{0.1}\right)^2\right)$. \\
The corresponding MFG system is,
\begin{equation}\label{1k}
\left\{
\begin{array}{rrrrr}
\phi_t+\nu \Delta \phi-\frac{1}{2}||\phi_x||^2+(1-\rho)\phi_x=0\\
 \rho_t -\nu \Delta \rho-\operatorname{div}((\phi_x-(1-\rho))\rho)=0\\
\rho(0, x)=0.05-0.9 \exp(\frac{-1}{2}(\frac{x-0.5}{0.1})^2), \\ \phi(T, x)=0. 
\end{array}
\right.
\vspace{-0.25cm}
\end{equation}\\

\begin{figure}
\centering\includegraphics[width=5.8cm]{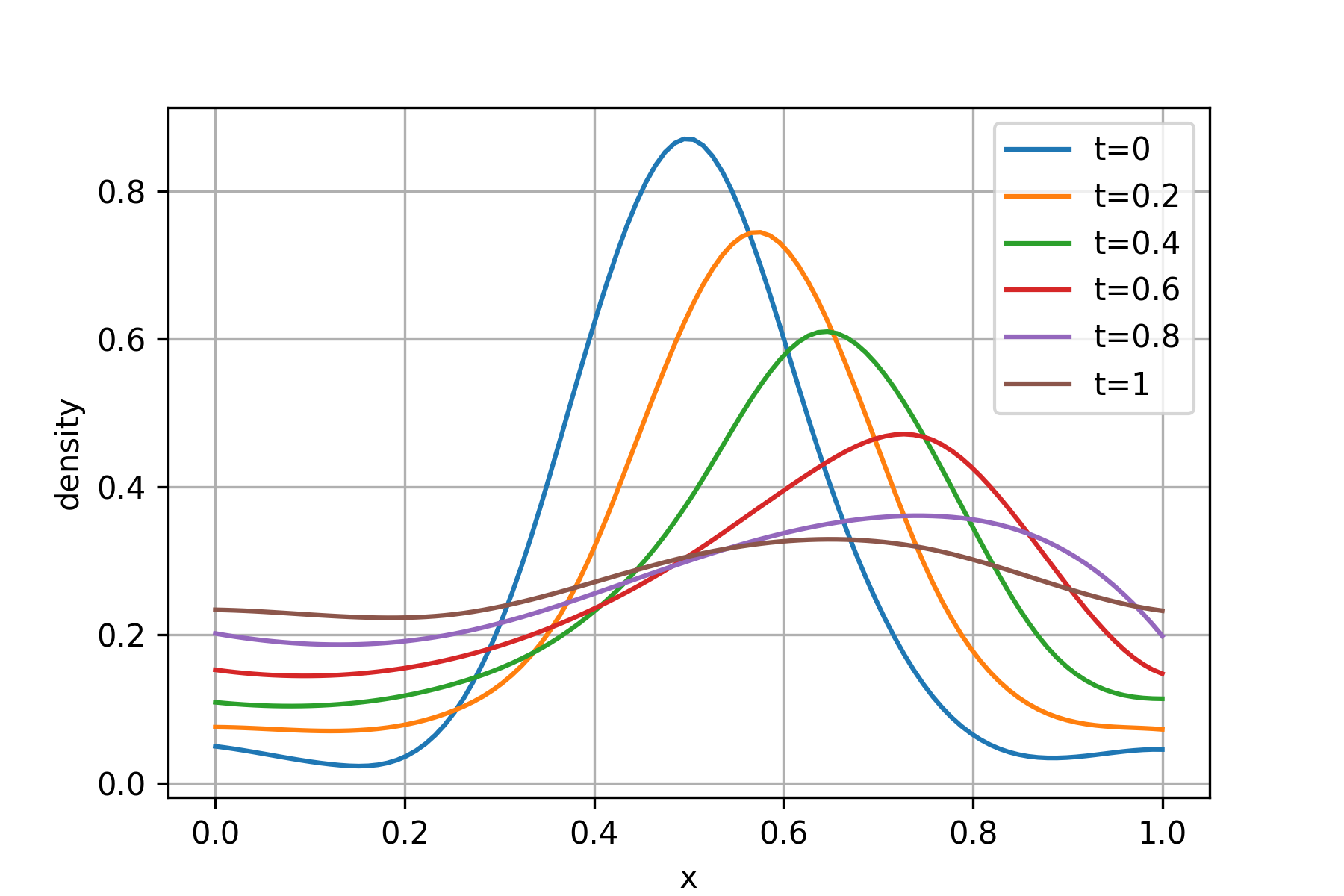}
\centering\includegraphics[width=5.8cm]{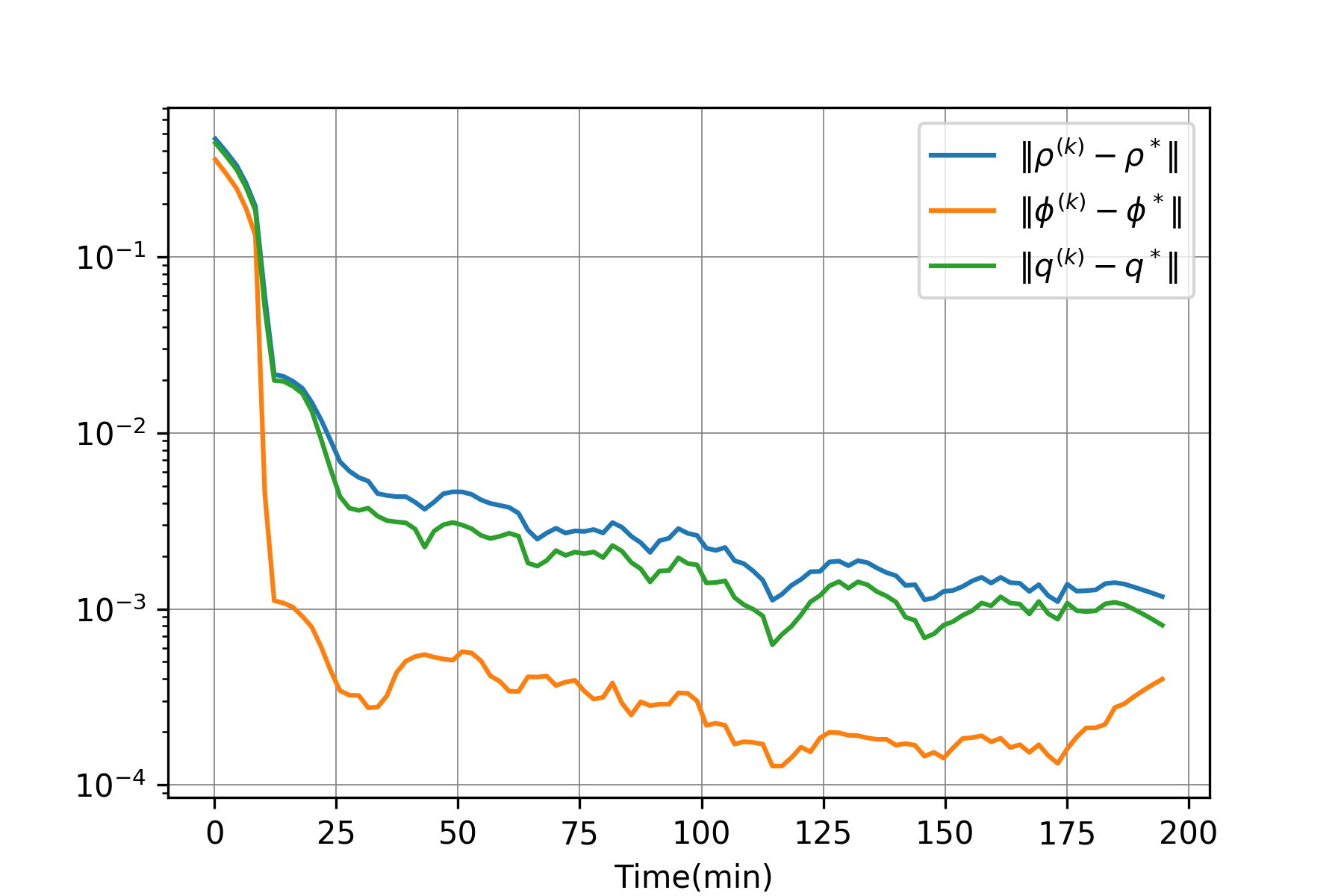}
\caption{We employ the DPI method to solve the problem (\ref{1k}). The figure shows the density solution at various time on the right and, on the left, the $L^{\infty}$ distance between $\rho^{(k)}$, $\phi^{(k)}$, and $q^{(k)}$ from DPI  and the final solution $\rho^{*}$, $\phi^{*}$, and $q^{*}$ obtained from the fixed-point algorithm.}
\label{fig222}
\end{figure}

{\bf Test 1:} In this test, we aim to solve (\ref{1k}) with periodic boundary conditions, we utilized Algorithm \ref{alg1} with a minibatch size of 50 samples at each iteration to obtain results. Our approach involved the use of neural networks with different hidden layers, each consisting of 100 neurons. Specifically, for $N_{\theta}$, we employed three hidden layers with the Gelu activation function, while for $N_{\omega}$ and $N_{\tau}$, we used a single hidden layer with the Sin activation function. During the training process, we employed the ADAM optimizer with a learning rate of $10^{-4}$ and weight decay of $10^{-3}$. To construct the neural networks, we adopted the ResNet architecture, incorporating a skip connection weight of 0.5.\\
To assess the performance of the methods, we conducted an analysis of the numerical results presented in Figure \ref{fig222}. This figure provides visual representations of the density solution as well as the $L^{\infty}$ Distance between $\rho^{(k)}$,  $\phi^{(k)}$  and $q^{(k)}$ from DPI and the final solution $\rho^{*}$,  $\phi^{*}$  and $q^{*}$ from fixed point algorithm. The  $L^{\infty}$  Distance values were obtained after applying a Savgol filter to enhance the clarity of the curves.
\\

\begin{figure}
\centering\includegraphics[width=5.8cm]{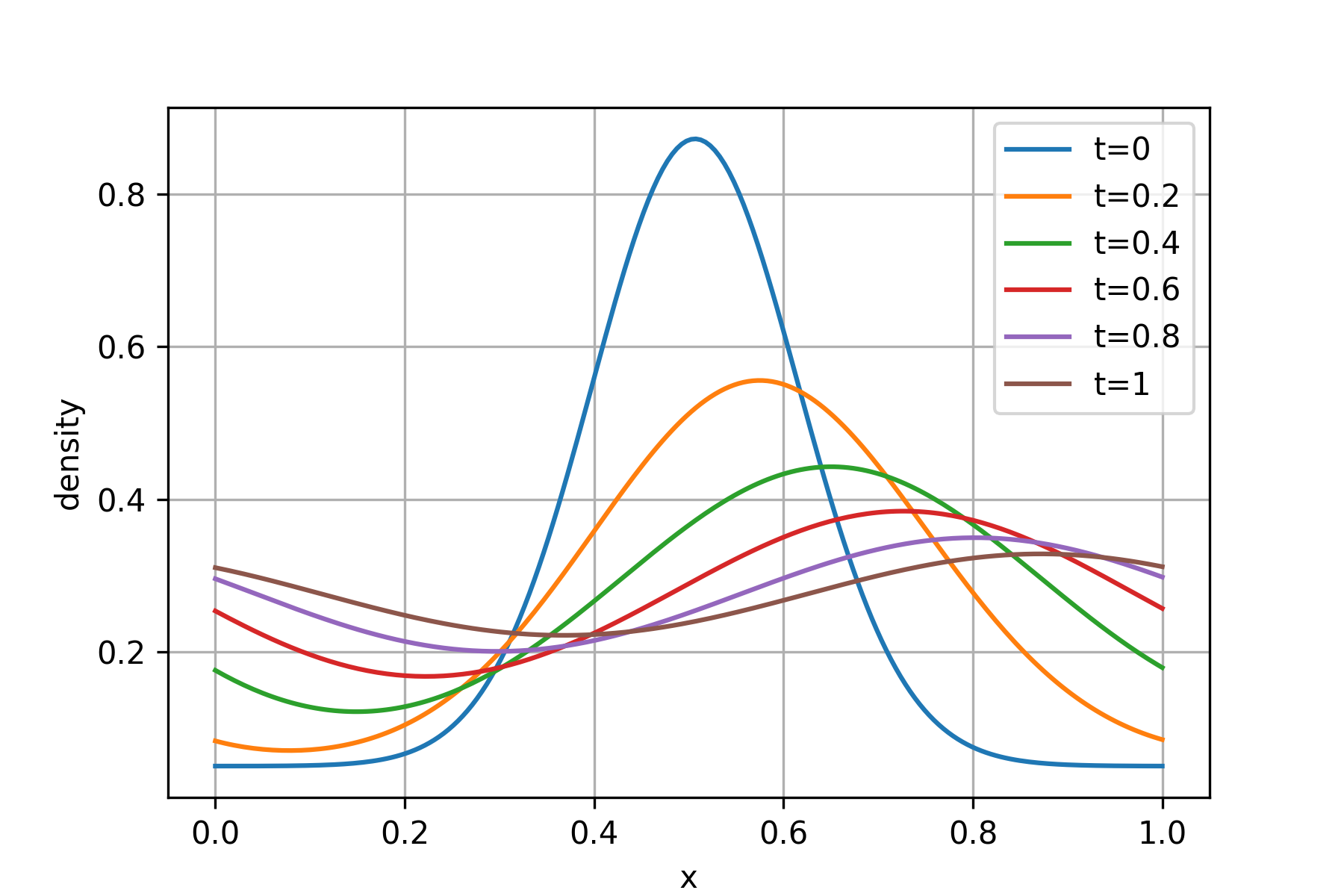}
\centering\includegraphics[width=5.8cm]{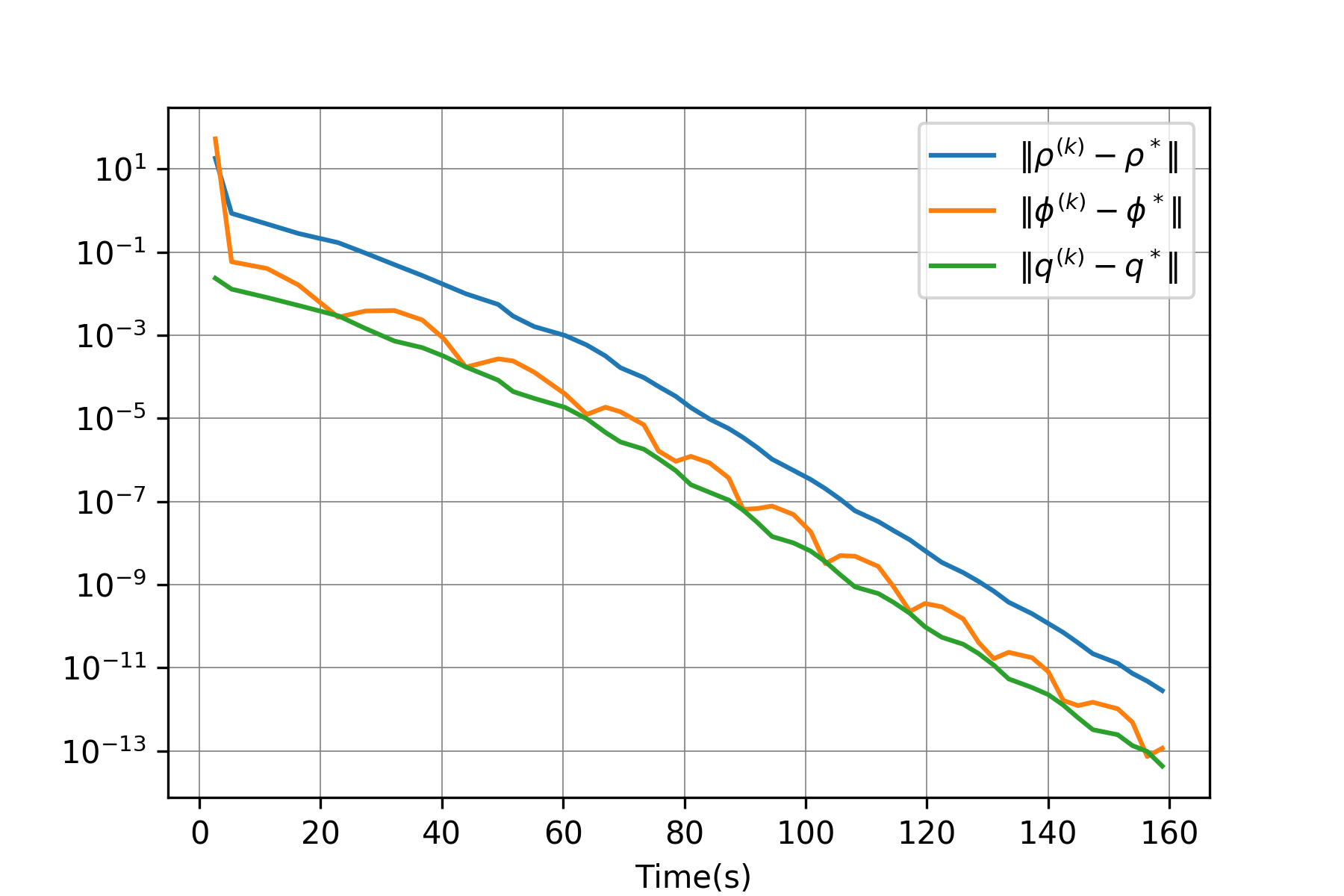}
\caption{We utilize the PI method to solve problem (\ref{1k}). The figure illustrates the density solution at various times on the right, while on the left, it demonstrates the $L^{\infty}$ distance between $\rho^{(k)}$, $\phi^{(k)}$, and $q^{(k)}$ from PI and the final solution $\rho^{}$, $\phi^{}$, and $q^{*}$ obtained from the fixed-point algorithm. The results indicate the effectiveness of PI for non-separable Hamiltonians in low dimensions} 
\label{figp222}
\end{figure}

{\bf Test 2:} The purpose of this test is to utilize the PI method in order to solve (\ref{1k}) with periodic boundary conditions. Similar to Test 3, we employ PI algorithm to solve the fully discretized (\ref{1k}) as follows: we start with an initial guess  $q^{(0)}_n: \mathcal{G} \to \mathbb{R}^{2d}$ for $n = 0, . . . , N,$ and initial and final data $\rho_0, \phi_N: \mathcal{G} \to \mathbb{R}$. We then iterate on $k \geq 0$ until convergence:\\
(i) Solve for $n=0, \ldots, N-1$ on $\mathcal{G}$
$$
\left\{\begin{array}{l}
\rho_{n+1}^{(k)}-d t\left( \Delta_{\sharp} \rho_{n+1}^{(k)}+\operatorname{div}_{\sharp}\left(\rho_{n+1}^{(k)} q_{n+1}^{(k)}\right)\right)=\rho_n^{(k)} \\
\rho_0^{(k)}=\rho_0
\end{array}\right.
$$
(ii) Solve for $n=N-1, \ldots, 0$ on $\mathcal{G}$
$$
\left\{\begin{aligned}
\phi_n^{(k)} & -d t\left( \Delta_{\sharp} \phi_n^{(k)}-q_{n, \pm}^{(k)} \cdot D_{\sharp} \phi_n^{(k)}\right) \\
& =\phi_{n+1}+\frac{dt}{2}\left(\left|q_{n+1, \pm}^{(k)}\right|^2+ \left|1 - \rho_{n+1}^{(k)}\right|^2 +\left(1 - \rho_{n+1}^{(k)}\right)q_{n+1, \pm}^{(k)}\right) \\
\phi_N^{(k)} & =\phi_N
\end{aligned}\right.
$$
(iii) Update the policy $q_n^{(k+1)}=D_{\sharp} \phi_n^{(k)}$ on $\mathcal{G}$ for $n=0, \ldots, N$, and set $k \leftarrow$ $k+1$.
In the following test, we choose  $\gamma=0$,  T = 1 for the final time, and $K = 50$. The grid consisted of $I=200$ nodes in space and $N=200$ nodes in time. The initial policy was initialized as $q^{(0)}_n \equiv (0, 0)$ on $\mathcal{G}$ for all $n$. Figure \ref{figp222} presents the numerical results, offering visual depictions of both the density solution and the $L^{\infty}$ Distance between $\rho^{(k)}$,  $\phi^{(k)}$  and $q^{(k)}$ from policy iteration and the final solution $\rho^{*}$,  $\phi^{*}$  and $q^{*}$ from fixed point algorithm.\\

The results of this study suggest that PI outperforms DPI in terms of effectiveness and speed. Consequently, it can be inferred that classical methods are often favoured for addressing such problems. However, due to the curse of dimensionality, deep learning approaches provide an alternative solution for high-dimensional cases.

\section{Conclusion}\label{Conclusion}
This paper presents Deep Policy Iteration (DPI), a novel approach designed to tackle computational challenges in high-dimensional stochastic Mean Field Games (MFG). DPI integrates the Deep Galerkin Method (DGM) with the stability and convergence benefits of the Policy Iteration Algorithm. Through the iterative training of three neural networks, DPI effectively approximates MFG solutions, overcoming the limitations posed by the curse of dimensionality in standard Policy Iteration. This extension enables the application of MFG solutions to high-dimensional scenarios, expanding the scope of this approach to MFG in the separable case, as demonstrated through numerical examples.

\appendix

 \bibliographystyle{elsarticle-num} 
 \bibliography{elsarticle-template-num}





\end{document}